# CLOSED FORM EXPRESSIONS FOR BAYESIAN SAMPLE SIZE

By B. Clarke and Ao Yuan

*University of British Columbia and Howard University*

Sample size criteria are often expressed in terms of the concentration of the posterior density, as controlled by some sort of error bound. Since this is done pre-experimentally, one can regard the posterior density as a function of the data. Thus, when a sample size criterion is formalized in terms of a functional of the posterior, its value is a random variable. Generally, such functionals have means under the true distribution.

We give asymptotic expressions for the expected value, under a fixed parameter, for certain types of functionals of the posterior density in a Bayesian analysis. The generality of our treatment permits us to choose functionals that encapsulate a variety of inference criteria and large ranges of error bounds. Consequently, we get simple inequalities which can be solved to give minimal sample sizes needed for various estimation goals. In several parametric examples, we verify that our asymptotic bounds give good approximations to the expected values of the functionals they approximate. Also, our numerical computations suggest our treatment gives reasonable results.

**1. Introduction.** Suppose $X^n = (X_1, \ldots, X_n)$ is IID $p(\cdot|\theta)$, where the $d$-dimensional parameter $\theta$ ranging over $\Theta \subset R^d$ is equipped with a prior probability $W(\cdot)$ having density $w(\theta)$ with respect to Lebesgue measure. Given an outcome $x^n = (x_1, \ldots, x_n)$ of $X^n$, Bayesian inference is based on the posterior density $w(\theta|x^n) = w(\theta)p(x^n|\theta)/m(x^n)$, where $m(x^n) = \int w(\theta)p(x^n|\theta)\,d\theta$ is the mixture density. Once a prior, likelihood and parametrization for $\theta$ are specified, the main pre-experimental task is to choose the sample size $n$. The size of $n$ will depend on the degree of accuracy desired and on the sense in which that accuracy is to be achieved.

Sample size determination in the Bayesian setting is an important and practical problem. As yet there is no general and accepted asymptotically valid closed form expression, such as we give here, that can be readily used to give minimally necessary sample sizes to achieve pre-specified inference objectives, even in seemingly simple cases. For instance, it has taken a series of papers (see [19] and the references therein) to provide a reasonable treatment for the difference of two proportions with independent Beta densities under a variety of criteria.

The lack of general expressions may be, in part, because the inferential criteria that have been used fall into three distinct classes. First, in the absence of a loss function, one often looks at properties of credibility sets— average length of the highest posterior density regions for instance. While this is often reasonable, the downside is that criteria that look for the worst case scenario often require overlarge sample sizes; see [14]. One way to correct for this is to include the cost of sampling in the optimality criterion.

Second, when a loss function is available, the decision theoretic approach originated by Raiffa and Schlaifer [20] can be used. One benefit of this approach is that it is easy to include the cost of sampling. The decision theoretic approach was developed in [18]. See also [1] and [16] for an information perspective; Pham-Gia and Turkkan ([19], Section 4) provided some general comments. Cheng, Su and Berry [3] established asymptotic expressions for sample size computation in the clinical trial context for dichotomous responses. A general discussion of the relative merits of decision theoretic approaches to sample size problems can be found in [14, 17, 18].

A third class of treatments of the sample size problem is more "evidentiary": These techniques tend to be based on hypothesis testing criteria such as Bayes factors (see [6, 7, 15]) or robustness; see [8]. The predictive probability criterion of [9], the distance between the posterior predictive density and the density updated on additional observations, and the direct evaluation of probabilities of events in the mixture distribution (see [4]) fall into this conceptual class as well. Since Bayesian testing can be framed as a decision problem, this third class can be regarded as a special case of the second class. However, the emphasis is different. Decision theoretic approaches tend to emphasize risks and expectations, while evidentiary approaches tend to focus on conditional probabilities, often posterior probabilities of hypotheses.

Because of this multiplicity of mathematically challenging criteria, it is not easy to parallel frequentist formulations. Nevertheless, many of these criteria can be represented as functionals $F$, not in general linear, of the posterior distribution $W(\cdot|X^n)$. For such cases, we provide a unified framework, indicating how it can be adapted to various settings.

Our overall goal is to give simple closed form asymptotic expressions in the form of inequalities that can be solved to give sample sizes. The reader interested primarily in these expressions can find four of them in Section 4, noted (APVC), (ACC), (ALC) and (ES), to indicate the criteria. [Expressions for similar cases are in Theorem 3.3 and in the Appendix; see (A.10), (A.11) and (A.13).] Informally, our central strategy for obtaining these expressions is the standard technique of approximating the leading term in an expansion of the expectation of a functional. Recall that $W(\cdot|X^n)$ is asymptotically $\Phi_{\hat{\theta},(nI(\theta))^{-1}}(\cdot)$ under $P_\theta$ in an $L^1$ sense. Here, $\Phi_{\mu,\Omega}(\cdot)$ is the distribution function for a Normal$(\mu,\Omega)$, with density denoted $\phi_{\mu,\Omega}(\cdot)$, and $\hat{\theta}$ is the maximum likelihood estimator (MLE), with asymptotic variance at a value $\theta$ given by the positive definite inverse Fisher information matrix $I(\theta)^{-1}$. If $\theta_0$ is the data generating parameter, adding and subtracting $E_{\theta_0}F(\Phi_{\hat{\theta},(nI(\theta_0))^{-1}}(\cdot))$ gives

$$(1.1) \qquad E_{\theta_0}F(W(\cdot|X^n)) = E_{\theta_0}F(\Phi_{\hat{\theta},(nI(\theta_0))^{-1}}(\cdot)) + E_{\theta_0}R_n(F),$$

where $R_n(F) = [E_{\theta_0}F(W(\cdot|X^n)) - E_{\theta_0}F(\Phi_{\hat{\theta},(nI(\theta_0))^{-1}}(\cdot))]$ is the remainder term and $F$ is a functional on distributions, that is, for any distribution $Q$, $F(Q) \in \mathbb{R}$. Our hope is that the remainder term will be small enough compared to the difference of the other two terms that (1.1) will permit asymptotically valid closed form expressions for the sample size criterion encapsulated by $F$.

1.1. *An example of the techniques.* Our verification that the remainder term in quantities like (1.1) is typically small rests on the foundational work of Johnson [10, 11], who developed Edgeworth style approximations for the posterior and certain posterior derived quantities such as percentiles and moments. Indeed, Edgeworth expansions and Johnson-style asymptotic expressions provide asymptotic control for the values of both terms on the right-hand side in (1.1), as $n \to \infty$, for various choices of $F$.

To see how these asymptotic expressions can be used to approximate the leading term of (1.1), and that the remainder term can be small compared to it, consider the following example. It is paradigmatic of our approach in its use of Johnson and Edgeworth expansions. The specific result can be obtained more readily by other techniques; however, our point is only to exemplify the reasoning informally.

Set $F(W(\cdot|X^n)) = F_\alpha(W(\cdot|X^n)) = W(D_n|X^n)$, where $D_n = (-\infty, a_n(\alpha))$ and $a_n = a_n(\alpha) = a_n(\alpha, X^n)$ is the $\alpha$th quantile under the posterior distribution $W(\cdot|X^n)$. Next, set

$$D'_n = \left(-\infty, \frac{1}{\sqrt{nI(\theta_0)}}\Phi^{-1}(\alpha) + Z_n\right] \equiv (-\infty, b_n],$$

in which $Z_n$ is an asymptotically standard normal random sequence of random variables. It is seen that $D'_n$ is the region corresponding to $D_n$ but under $\Phi_{Z_n,(nI(\theta_0))^{-1}}(\cdot)$, in which we have used $Z_n$ in place of $\hat\theta$ by asymptotic normality of the MLE. That is, $D'_n$ approximates $D_n$. In this case, the first term on the right-hand side of (1.1) is

$$
\begin{aligned}
&E_{\theta_0}\Phi_{Z_n,(nI(\theta_0))^{-1}}(D'_n)\\
(1.2)\qquad &= E_{\theta_0}\left(\int_{-\infty}^{\Phi^{-1}(\alpha)/\sqrt{nI(\theta_0)}+Z_n}\frac{\sqrt{n}I^{1/2}(\theta_0)}{\sqrt{2\pi}}e^{-(1/2n)I(\theta_0)(\theta-Z_n)^2}\,d\theta\right)\\
&= \frac{1}{\sqrt{2\pi}}\int_{-\infty}^{\Phi^{-1}(\alpha)}e^{-t^2/2}\,dt = \alpha.
\end{aligned}
$$

The remainder term in (1.1) is

$$(1.3)\qquad E_{\theta_0}R_n = E_{\theta_0}\chi_{(a_n\wedge b_n, a_n\vee b_n)}(\cdot) = E_{\theta_0}|a_n - b_n|.$$

Posterior normality suggests $(1.3)\to 0$, but we want a rate that is small relative to the rate of convergence of the left-hand side of (1.1) to (1.2) which we take to be $o(1)$. We ignore details on this latter rate since it is not the point. Now, to get a rate for $(1.3)\to 0$, we use a modification of Johnson ([11], Theorem 5.1); it is justified below in Theorem 2.1. Thus, we have that quantiles such as $a_n$ satisfy

$$a_n = (nI(\theta_0))^{-1/2}\left[\Phi^{-1}(\alpha) + \sum_{j=1}^J \tau_j(\alpha)n^{-j/2} + O(n^{-(J+1)/2})\right] + \hat\theta_n,$$

where the $\tau_j$'s are polynomials with bounded coefficients that depend on the data $X^n$, and $J\geq 1$. Now, we can write

$$
\begin{aligned}
E_{\theta_0}|a_n - b_n| &= E_{\theta_0}\Bigg|n^{-1/2}I^{-1/2}(\theta_0)\bigg[\Phi^{-1}(\alpha) + \sum_{j=1}^J \tau_j(\alpha)n^{-j/2}\\
&\qquad\qquad\qquad\qquad\qquad + O(n^{-(J+1)/2})\bigg] + \hat\theta_n\\
&\qquad\qquad - (n^{-1/2}I^{-1/2}(\theta_0)\Phi^{-1}(\alpha) + Z_n)\Bigg|\\
(1.4)\quad &= E_{\theta_0}|\hat\theta_n - Z_n| + O(n^{-1/2})\\
&\leq E_{\theta_0}|\hat\theta_n - \theta_0| + E_{\theta_0}|Z_n - \theta_0| + O(n^{-1/2})\\
&= n^{-1/2}I^{-1/2}(\theta_0)(E_{\theta_0}|\sqrt{n}I^{1/2}(\theta_0)(\hat\theta_n - \theta_0)|\\
&\qquad\qquad\qquad + E_{\theta_0}|\sqrt{n}I^{1/2}(\theta_0)(Z_n - \theta_0)|)\\
&\quad + O(n^{-1/2}).
\end{aligned}
$$

Expression (1.4) can be controlled by using an Edgeworth expansion for the density of $\hat{\theta}$ under $\theta_0$ in the first term in parentheses, namely, $E_{\theta_0}\sqrt{n}I^{1/2}(\theta_0)(\hat{\theta}_n - \theta_0)$. Using this approximation and recognizing limiting normal forms gives that, term by term, (1.4) is

$$n^{-1/2}I^{-1/2}(\theta_0)\left(\int |z|\phi(z)\,dz + \sum_{k=1}^{K} n^{-k/2}\int |z|P_k(z)\,dz \right.$$
$$\left. + o(n^{-K/2})\int \frac{|z|}{1+|z|^{K+2}}\,dz + \int |z|\phi(z)\,dz \right) + O(n^{-1/2}).$$

So, (1.3) is $O(1/\sqrt{n})$ and the left-hand side of (1.1) is

(1.5) $\qquad E_{\theta_0}F_\alpha(W(\cdot|X^n)) = \alpha + o(1) + O(n^{-1/2}),$

that is, the expected Bayesian coverage probability is always $\alpha + o(1)$.

Improving (1.5) leads to inequalities that can be solved to give sample sizes. That is, careful use of the Edgeworth and Johnson expansions that we used to control (1.3) and (1.4) will give an error term of order $o(1/\sqrt{n})$. So, we can find $N = N(\varepsilon)$ large enough that, for a specified range of parameter values $\theta$, we would have $|E_\theta F_\alpha(W(\cdot|X^n)) - \alpha| < \varepsilon$ for $n > N$. Details on this case are given below in Example 3 of Section 4. The "nicest" cases occur when the first term in (1.1) is independent of the value of $\hat{\theta}$ and the second term goes to zero. As suggested by the form of (1.2), when the first term in (1.1) depends on an estimator such as $a_n$ or $\hat{\theta}$, we expect an asymptotically normal random variable $Z_n$ to appear in the limit. In these cases, we want the second term of (1.1) to go to zero at a fast enough rate. Thus, we want to give an expansion for it as a sum of powers of $1/\sqrt{n}$ times evaluations of expectations.

1.2. *Expected values of functionals of the posterior.* Before proceeding with the mathematical formalities, we suggest that the formulation we have adopted here—representing sample size criteria as expectations of functionals of the posterior—is the right one, in the sense that it is general enough to encapsulate all the important cases, yet narrow enough to permit straightforward analysis and use.

The three classes identified earlier—Bayes credibility, decision theoretic and evidentiary—suggest that many authors have, implicitly or explicitly, studied criteria that amount to functionals of the posterior, if not expectations of them. Indeed, the pure Bayes and evidentiary approaches amount to studying functionals of the posterior and most of the decision theoretic optimality criteria can be written as functionals of the posterior; most often these are clearly expectations. Moreover, taking expectations over the sample space pre-experimentally is standard Bayesian practice for design problems. This is done in [23], for instance, an approach that motivated the present work. Wang and Gelfand proposed a simulation based technique for determining a sample size large enough to achieve various pre-experimentally specified criteria.

All the criteria used in [23] are special cases of the form $E(T(Y)) \leq \varepsilon$, where $T$ is a nonnegative function in which the data $Y$ appears via conditioning; see [23], Section 2, equation (6). Their simulation technique has a broad scope of application, and should be at least as accurate as approximations based on asymptotic expansions. The special cases of $F$ we use here are taken from [23].

We comment that some of the criteria used in Wang and Gelfand's simulations, for instance, the average cover criterion, ACC, and average length criterion, ALC, have been studied mathematically. For instance, Joseph and Bélisle [12] and Joseph, du Berger and Bélisle [13] derived inequalities the sample size must satisfy under certain prior specifications for normal and binomial models. Wang and Gelfand's work [23] is important because these special cases may not cover all the settings of interest.

Unfortunately, simulations may not always be easy to do. Moreover, the distinction between the sampling and fitting priors used in [23] may be a layer of conservatism that is not necessary. Aside from computational ease, Sahu and Smith ([21], Section 2.3) argue that using sampling and fitting priors permits weaker assumptions for the validity of inference. However, one could use a single objective prior for both sampling and fitting purposes to achieve essentially the same inferential validity. In either case, there remains a role in Bayesian experimental design for a good closed form expression for sample sizes.

Expression (1.1) suggests a different tack for obtaining the kind of closed form expressions we want. One could approximate $E_\theta F(W(\cdot|X^n))$ by $E_\theta F(\hat{N}(\hat{\theta}, (n\hat{I}(\hat{\theta}))^{-1}))$, where $\hat{N}$ is a Laplace approximation to the posterior, instead of a Johnson style expansion. The two approaches—Johnson and Laplace—probably require similar hypotheses. Arguably, the Laplace expansion is conceptually easier. However, Johnson expansions give an approximation to $F(W(\cdot|X^n))$ directly rather than separately approximating $F$ and $W(\cdot|X^n)$. One could use more terms in the Laplace approximation, evaluate $F$ on those terms, and then approximate $F$, but the complexity would likely exceed what we have done here. The Johnson expansions are readily available and more direct, although a confirmatory treatment using Laplace's method would be welcome.

The structure of this paper is as follows. Section 2 gives the theoretical context of our work: We observe generalizations of key results in Johnson [11] and state the version of Edgeworth expansions we will need. Then, we give a simple result, Proposition 2.1, that formalizes the strategy implicit in (1.1). It seems that getting an asymptotic expression for general functionals $F$ is a hard problem so, in Section 3, we give asymptotic expressions for three kinds of terms that often arise in special cases of functionals of the posterior density. Two of these theorems are derived from [11], and one is new. The most technical arguments from this section are relegated to the Appendix at the end. Section 4 uses our main results to show how four established criteria for sample size determination admit asymptotically valid closed form expressions. In Section 5 we compare the results of our asymptotic expressions to closed form expressions obtained from three exponential families equipped with conjugate priors. It is seen that our asymptotic expansions typically match the leading $1/\sqrt{n}$ terms in those cases. In addition, Section 5 presents numerical results which confirm our approximations are reasonably accurate.

**2. Theoretical context.** We consider the case that $F$ is a functional on distributions such as the posterior $W(\cdot|X^n = x^n)$ for a parameter. We assume $F$ represents something about how distributions concentrate at a specific value in their support. Our interest here focuses on the class of $F$ only in that we want to include the commonly occurring sample size criteria used in [23].

We will need two assumptions to control the leading term in an expansion for $E(F)$. The first is drawn from [11], Theorem 2.1: The expectation of the functional of the posterior, $EF(W(\cdot|X^n))$ minus its normal approximation [see (1.1)] must have an expansion of the form established by Johnson [11]. The second assumption is that the classical Edgeworth expansion can be used to approximate the sampling distribution of $\hat{\theta}_n$ when $\theta$ is taken as true.

To begin, we make Assumptions 1–9 in [11], modifying them only by permitting $\theta$ to range over a set $\Omega \subset \mathbb{R}^d$. Together, these are the standard "expected local sup" conditions that ensure the consistency, asymptotic normality and efficiency of the MLE. Assumption 8, for instance, bounds the first two derivatives of $\log p(x|\theta)$ by an integrable function so that, when $d = 1$,

$$\hat{I}(\hat{\theta}) = -\frac{1}{n}\sum_{i=1}^{n}\frac{\partial^2}{\partial\theta^2}\log p(X_i|\hat{\theta}) \stackrel{\text{a.s.}}{\to} -E_{\theta_0}\frac{\partial^2}{\partial\theta^2}\log p(X|\theta) = I(\theta),$$

which generalizes directly to multivariate $\theta$.

To set up our first result, we need some notation. Let $\theta$ be a random realization of $\Theta$, $\hat{\phi}_n = \sqrt{n}\hat{I}^{1/2}(\hat{\theta}_n)(\theta - \hat{\theta}_n)$ and consider Johnson expanding the posterior distribution function $\hat{W}(\hat{\phi}|X^n)$ of $\hat{\phi}_n$. Johnson [11] obtained an expansion for $\hat{W}(\hat{\phi}_n|X^n)$ in terms of normal densities with polynomial factors when $\theta$ is one-dimensional. The expansion uses $(n\hat{I}(\hat{\theta}_n))^{-1}$ as the empirical variance of $\hat{\theta} - \theta$ and holds in an almost sure sense, for $n > N_x$, where $N_x$ depends on the observed sample $x = x^n$. This is almost the expansion we want. For our purpose, we set $\psi = \psi_n = \sqrt{n} I^{1/2}(\theta_0)(\theta - \hat{\theta}_n)$ for given $\hat{\theta}_n$ and denote the posterior distribution function of it by $W_o(\cdot | X^n)$. Writing the distribution of the $d$-dimensional standard normal $N(\mathbf{0}, I_d)$ as $\Phi(\cdot)$, with density $\phi(\cdot)$, we have $\Phi(\sqrt{n} I^{1/2}(\theta_0)(\theta - \hat{\theta}_n)) = \Phi_{\hat{\theta}_n, I^{-1}(\theta_0)/n}(\theta)$ and $\phi(\sqrt{n} I^{1/2}(\theta_0)(\theta - \hat{\theta}_n)) = |nI(\theta_0)|^{-1/2} \phi_{\hat{\theta}_n, I^{-1}(\theta_0)/n}(\theta)$. Let $w^{(r)}(\theta)$ be the $r$th (vector) derivative of the prior density $w(\theta)$, when it exists, and write $\hat{I}_r(\theta) = \frac{1}{n|r|!} \sum_{i=1}^n \frac{\partial^{|r|}}{\partial \theta^r} \log p(X_i | \hat{\theta})$ for a vector $r = (r_1, \ldots, r_d)$, where $|r| = k$ means $r_1 + \cdots + r_d = k$, and for $\theta = (\theta_1, \ldots, \theta_d)$, $\theta^r$ means $\theta_1^{r_1} \cdots \theta_d^{r_d}$. Examination of [11] gives the following.

THEOREM 2.1.  *Suppose all derivatives of $\log p(\cdot | \theta)$ of order $J + 3$ or less exist and are continuous and that all the derivatives $|(\partial^{|r|}/\partial \theta^r) \log p(x|\theta)|$, for $|r| \leq J + 3$, are bounded in an open set containing $\theta_0$ by a function $G(x)$ with $EG(X)$ finite. Suppose also that all derivatives of $w$ up to order $J + 1$ exist and are continuous in a neighborhood of $\theta_0$. Then, for given $\theta_0$, there are a sequence of sets $S_n$ with $P_{\theta_0}(S_n^c) = o(1)$, and an integer $N$, so that, for $x^n \in S_n$, Theorems 2.1, 3.1, 4.1, 5.1 and 5.2 of [11] continue to hold with $\hat{W}(\phi|X^n)$ replaced by $W_o(\psi|X^n)$ when $n \geq N$. That is, we have:*

(A) *For the posterior distribution:*

(2.1)
$$\left| W_o(\psi|X^n) - \Phi(\psi) - \sum_{j=1}^J n^{-j/2} \phi(\psi) \gamma_j(\psi, X^n) \right| \leq C n^{-(J+1)/2},$$

$$n > N, X^n \in S_n,$$

*where $C > 0$ is a constant, and the $\gamma_j(\psi)$'s are polynomials in $\psi$ with bounded coefficients.*

(B) *For posterior moments: For each integer $i \leq K - 1$, there are a sequence of functions $\{\lambda_{ij}(X^n)\}$, a constant $C > 0$ and an integer $N_i$ so that*

(2.2)
$$\left| E_{W_o(\cdot|X^n)}(I_{S_n} I^{i/2}(\theta_0)(\theta - \hat{\theta}_n)^i) - \sum_{j=i}^J \lambda_{ij}(X^n) n^{-j/2} \right| \leq C n^{-(J+1)/2},$$

$$n > N_i,$$

*on a set $S_n(i)$ with $P_{\theta_0}(S_n(i)^c) \to 0$, where $\lambda_{ij}(X^n) = 0$ for $j$ odd, and for $i$ even we have*

$$\lambda_{ii}(X^n) = 2^{i/2} \Gamma((i+1)/2)/\Gamma(1/2),$$

*while for $i$ odd we have*

$$\lambda_{i,i+1}(X^n) = 2^{(i+1)/2}(2(i+1)I_{3n}(\hat{\theta}_n)\Gamma((i+4)/2) \\ + \Gamma((i+2)/2)w^{(1)}(\hat{\theta}_n)/w(\hat{\theta}_n))/\Gamma(1/2),$$

*all of which are bounded in $X^n$.*

(C) *For inverse quantiles: Let $\eta(\xi) = \Phi^{-1}(W_o(\xi|X^n))$ be the transformed quantile of $W_o(\cdot|X^n)$. Then*

$$(2.3) \quad \left|\eta(\xi) - \xi - \sum_{j=1}^{J} n^{-j/2}\omega_j(\xi)\right| \leq Cn^{-(J+1)/2}, \qquad n > N, X^n \in S_n,$$

*where $C > 0$ is a constant, for some functions $\omega_j(\xi) = \omega_j(\xi, X^n)$ that are polynomials in $\xi$ with coefficients bounded for large enough $n$.*

(D) *For posterior quantiles: For a solution $\eta = \Phi^{-1}(W_o(\xi(\eta)|X^n))$, we have the following:*

(i)

$$(2.4) \quad \left|\xi_n(\eta) - \eta - \sum_{j=1}^{J} n^{-j/2}\tau_j(\eta)\right| \leq Cn^{-(J+1)/2}, \qquad n > N, X^n \in S_n,$$

*where $C > 0$ is a constant and the functions $\tau_j(\cdot)$ are polynomials in $\eta$ with bounded coefficients.*

(ii) *If we set $\eta = \alpha$th percentile of $\Phi$, then*

$$(2.5) \left|W_o\left(\eta + \sum_{j=1}^{J} n^{-j/2}\tau_j(\eta)|X^n\right) - \alpha\right| \leq Cn^{-(J+1)/2}, \qquad n > N, X^n \in S_n.$$

REMARK. This collection of statements differs from Johnson's [11] results because we observe it for general $d$-dimensional parameters, a single choice of $N$ independent of the data string, and have replaced the empirical Fisher information by its population value in the standardization of the MLE. Replacing the $N_{k,x}$'s in [11] by a single fixed $N$ means we can only get a Johnson expansion valid for $x^n$ in a set $S_n$ with probability increasing as $P_{\theta_0}(S_n) = 1 - o(1)$. To ensure $P_{\theta_o}(S_n^c) = o(1)$, we will typically need laws of large numbers to hold for the $\hat{I}_r$'s occurring in the expansion; we assume these as needed. Faster rates for $P_{\theta_o}(S_n^c) \to 0$, for instance, $P_{\theta_o}(S_n^c) \leq e^{-n\gamma}$ for $\gamma > 0$, can be obtained by imposing moment generating function assumptions to get a large deviations principle.

Note that $I(\theta_0)$ is used in the standardization of the MLE, but the coefficients in the expansion remain empirical. That is, the coefficients in the polynomials of the expansions are functions of the data, usually estimates of population quantities of the form [11], equations (2.25) and (2.26). When it is important to replace these with differentiable quantities, as in the proof of Theorem 3.3, we will use approximations such as $\hat{I}(\hat{\theta}) = I(\theta_0) + o_p(1)$; the $o_p(1)$ term in such approximations is what limits the accuracy of our expansions.

PROOF OF THEOREM 2.1. Proofs for (2.1)–(2.5) are all modifications of the techniques in [11]. To demonstrate the modifications, consider (2.1). It will be enough to check the proof of Theorem 2.1 in [11] line by line.

First, the main difference due to the dimensionality is that occurrences of powers $(\theta - \hat{\theta}_n)^r$ in the one-dimensional case must be replaced by the multi-dimensional version, $\sum_{|r|=k}(\theta - \hat{\theta}_n)^r$ for a $d$-tuple nonnegative integer vector $r$.

Johnson used bounds $N_{k,x}, k = 1, \ldots, 5$, in his proof. The first two, $N_{1,x}$ and $N_{2,x}$, are used in his Lemmas 2.1 and 2.2, which are not needed in our case, since we are replacing $\hat{I}(\hat{\theta}_n)$ by $I(\theta_0)$ (Note that in the statement of Lemma 2.2 in [11], $f(x_i, \theta)$ in the denominator should be $f(x_i, \hat{\theta}_n)$.) The next two, $N_{3,x}$ and $N_{4,x}$, are from Lemmas 2.3 and 2.4. They arise from using the strong law of large numbers finitely many times to get inequalities. Denote the set on which the strong laws fail for a given $n$ by $S_n^c$. Then, the conclusions in Lemmas 2.3 and 2.4 hold for all $x \in S_n$, and $P(S_n^c) = o(1)$. This property of the strong law holds even when $\hat{I}(\hat{\theta})$ is replaced by $I(\theta_0)$. Finally, $N_{5,x} > N_{4,x}$ is used to allow the finite term approximations (2.21) and (2.22) to be used in the expansions (2.19) and (2.20). The sets of $x^n$'s on which this fails have probability tending to zero. Thus, they can be put into $S_n^c$ too, and $N$ can be chosen independent of $x^n$. □

It is seen from (2.1) that, for $n > N$ and $X^n \in S_n$,

$$W_o(\psi|X^n) = \Phi(\psi) + \sum_{j=1}^{J} n^{-j/2} \phi(\psi) \gamma_j(\psi, X^n) + n^{-(J+1)/2} \gamma_{J+1}(\psi, X^n),$$

for positive integers $J$, where the polynomials $\gamma_j(\psi)$ in $\psi$ have finite coefficients. Note that $\gamma_{J+1}$ is not known to be of the form of the $\gamma_j$'s when $j \leq J$; it is only known to be bounded. The other expansions (2.2)–(2.4) give analogous statements.

We formalize this class of posterior approximations in the following definition. First, we say that $P_W(x^n)$ is a posterior derived object if and only of $P_W(x^n)$ is a function of the posterior distribution $W(\cdot|x^n)$. Here, we have chosen $W_o(\cdot|X^n)$ as the form of the posterior for our work. The class of $P_W(x^n)$ does not matter, but the use of $W(\cdot|x^n)$ does. We rule out the appearance of parameters or their estimates apart from $I(\theta_0)$. Thus, the posterior itself and a posterior quantile are both posterior derived objects.

ASSUMPTION JE. A posterior derived object $P_{W_o}(x^n)$ is Johnson expandable of order $J$ if and only if it has a Johnson expansion of the following form: There are an $N$ and an $S_n$ with $P_{\theta_0}(S_n^c) = o(1)$ so that, for $n > N$, we have

$$\left| P_{W_o}(x^n) - \sum_{j=0}^{J} \frac{\gamma_j(x^n)}{n^{j/2}} \right| \leq \frac{C}{n^{(J+1)/2}},$$

for some $C > 0$, where the $\gamma_j(x^n)$'s are any quantities that depend only on $W_o(\cdot|x^n)$.

We assume that all Assumption JE's are nontrivial, that is, the $j = 0$ term is not $P_{W_o}(x^n)$.

Next, we turn to the other asymptotic expansion assumption we will need. For the MLE $\hat{\theta}_n$ of $\theta$ based on $p(X^n|\theta)$, let $f_n(\cdot) = f_n(\hat{\theta}|\theta)$ be the density function of $\hat{\theta}_n$ when $\theta$ is the true value, and let $g_n(\cdot) = g_n(\cdot|\theta)$ be the density of $T = T_n = \sqrt{n}I^{1/2}(\theta_0)(\hat{\theta}_n - \theta)$ given $\theta$. (It is seen that $T$ is a function of $\hat{\theta}$ for fixed $\theta$, whereas $\alpha_n$ is a function of $\theta$ for given $\hat{\theta}$.) Observe that

$$f_n(\theta) = |nI(\theta_0)|^{1/2} g_n(\sqrt{n} I^{1/2}(\theta_0)(\theta - \theta_0)).$$

So, to get an expansion for $f_n$, it is enough to get one for $g_n$. For later use, we record

$$\Phi_{\hat{\theta}_n, I^{-1}(\theta_0)/n}(\theta) = \Phi(\sqrt{n} I^{1/2}(\theta_0)(\theta - \hat{\theta}_n))$$

and

$$\phi_{\theta_0, (nI(\theta_0))^{-1}}(\theta) = |nI(\theta_0)|^{1/2} \phi_d(\sqrt{n} I^{1/2}(\theta_0)(\theta - \theta_0)).$$

The expansion for $g_n$ will depend on the form of the MLE. For many parametric families, $\hat{\theta}_n$ can be expressed as

$$\hat{\theta}_n = s\left( \frac{1}{n} \sum_{i=1}^{n} h(X_i) \right),$$

for some $s(\cdot)$ and $h(\cdot)$. Thus, as argued in [24], we often have

$$g_n(t) = \phi_d(t) + \sum_{k=1}^{K} n^{-k/2} P_k(t) + o(n^{-K/2}) \frac{1}{1 + \|t\|^{K+2}},$$

where the error $o(n^{-K/2})$ is uniform over $\theta$ in a compact set and $t = \sqrt{n} I^{1/2}(\theta) \times (\hat{\theta}_n - \theta)$. The $P_k(v)$'s are polynomials given by

$$\phi_d^{-1}(v) \sum_{q=1}^{k} \frac{1}{q!} \sum_{l_1 + \cdots + l_q = k, |r_m| = l_m + 2, (1 \leq m \leq q)} \sum \frac{\chi_{r_1} \cdots \chi_{r_q}}{r_1! \cdots r_q!}$$

$$\times (-1)^{|r_1| + \cdots + |r_q|} D^{r_1 + \cdots + r_q} \phi_d(v)$$

in which $\chi_r$, for a vector $r$, is the $r$th cumulant; see [2].

ASSUMPTION EE. The Edgeworth expansion of order $K$ for $f_n(\cdot)$ induced from $g_n(\cdot)$ is

$$f_n(\theta) = \phi_{\theta_0,(nI(\theta_0))^{-1}}(\theta) + \sum_{k=1}^{K} n^{-k/2} P_k(\sqrt{n}I^{1/2}(\theta_0)(\theta - \theta_0))\phi_{\theta_0,(nI(\theta_0))^{-1}}(\theta)$$

$$+ o(n^{-K/2}) \frac{|nI(\theta_0)|^{1/2}}{1 + \|\sqrt{n}I^{1/2}(\theta_0)(\theta - \theta_0)\|^{K+2}},$$

when it exists, where $\theta$ is a dummy variable varying over values of $\hat{\theta}$ and the error $o(n^{-(K-2)/2})$ is uniform for $\theta$ in a compact set.

We comment that Yuan and Clarke [24] do not prove Assumption EE in full generality. They only establish uniformity for the density of the mean and for a certain restricted class of functions of the mean. However, the discussion in [24] suggests that Assumption EE holds in much greater generality even though a formal proof does not yet exist. Indeed, when it fails, it seems to do so only on sets of very small probability which are enough to prevent the supremum from going to zero. Consequently, we suggest Assumption EE is an acceptable hypothesis in a design setting where we are primarily interested in average behavior rather than worst case behavior.

Note that Assumption EE permits us to take expectations over the parameter space and the sample space because the approximation is uniformly good over both $\theta$ and $X^n$. Indeed, Assumption EE immediately gives an expression for the mean of $\hat{\theta}$ because

$$\int \theta f_n(\theta) \, d\theta = \int |nI(\theta_0)|^{1/2} \theta \phi_d(\sqrt{n}I^{1/2}(\theta_0)(\theta - \theta_0)) \, d\theta$$

$$(2.6) \qquad + \sum_{k=1}^{K} n^{-k/2} \int |nI(\theta_0)|^{1/2} \theta P_k(\sqrt{n}I^{1/2}(\theta_0)(\theta - \theta_0))$$

$$\times \phi_d(\sqrt{n}I^{1/2}(\theta_0)(\theta - \theta_0)) \, d\theta$$

$$+ o(n^{-K/2}) \int \frac{|nI(\theta_0)|^{1/2}\theta}{1 + \|\sqrt{n}I^{1/2}(\theta_0)(\theta - \theta_0)\|^K} \, d\theta$$

$$= \theta_0 + \int u \phi_d(u) \, du$$

$$+ \sum_{k=1}^{J} n^{-k/2} \int (\theta_0 + u/(\sqrt{n}|I(\theta_0)|^{1/2})) P_k(u) \phi_d(u) \, du$$

$$+ o(n^{-K/2}) \int \frac{\theta_0 + u}{1 + \|u\|^K} \, du.$$

$$= \theta_0 + \sum_{k=1}^{J} n^{-k/2} \theta_0 P_k(\sigma)$$

$$+ \sum_{k=1}^{J} n^{-(k+1)/2} |I(\theta_0)|^{-1/2} P_{1,k}(\sigma) + o(n^{-K/2}),$$

where $P_k(\sigma)$ and $P_{1,k}(\sigma)$ are the expectations of $P_k(u)$ and $uP_k(u)$. The argument $\sigma$ signifies that powers $u^m$ are replaced by $\sigma_m$'s, the $m$th moments of $N(0,1)$. To see this, suppose $Z = (Z_1, \ldots, Z_d) \sim N(\mathbf{0}, I_d)$ and that the $i$th term in $P_k(u)$ has the form $a_i u_1^{i_1} \cdots u_d^{i_d}$. Then the term in its expectation is $a_i \int u(u_1^{i_1} \cdots u_d^{i_d}) \phi_d(u) \, du$, which equals $a_i E(Z_1^{i_1+1} Z_2^{i_2} \cdots Z_d^{i_d}, \ldots, Z_1^{i_1} \cdots Z_d^{i_d+1}) = a_i(\sigma_{i_1+1} \sigma_{i_2} \cdots \sigma_{i_d}, \ldots, \sigma_{i_1} \cdots \sigma_{i_d+1})$, a vector with entries in which the powers of $u_i$ correspond to standard normal moments.

Recall, our goal is to derive asymptotically, for pre-specified $\varepsilon > 0$ and $F$, the minimal sample size $n$ to achieve

(2.7) $$E_{\theta_0} F(W(\cdot|X^n)) \leq \varepsilon,$$

where the expectation is with respect to the density $p(x^n|\theta_0)$. Our main approach to (2.7) rests on the following general procedure for the computation of the asymptotic expected behavior of functionals of the posterior distribution. As indicated in the Introduction, let

(2.8) $$R_n = F(W(\cdot|X^n)) - F(\Phi_{\hat{\theta}_n, (nI(\theta_0))^{-1}}(\cdot)),$$

where, under $\theta_0$, $\hat{\theta}_n$ is distributed as in Assumption EE, and we have done the standardization in the limiting normal rather than in the nonstandardized posterior $W(\cdot|X^n)$ for $\theta$.

PROPOSITION 2.1. *Functionals of the posterior distribution function* $W(\cdot|X^n)$ *satisfy the following:*

(i) *If* $F(\Phi_{z,(nI(\theta_0))^{-1}}(\theta))$ *is independent of* $z$, *then if Assumption JE holds for some* $J \geq 1$, *we have*

(2.9) $$E_{\theta_0} F(W(\theta|X^n)) = F(\Phi_{\mathbf{0},(nI(\theta_0))^{-1}}(\theta)) + E_{\theta_0} R_n.$$

(ii) *If Assumption EE holds for some* $K \geq 1$, *we have that*

(2.10) $$E_{\theta_0} F(W(\theta|X^n)) = E_{\theta_0} F(\Phi(Z + \sqrt{n} I^{1/2}(\theta_0)(\theta - \theta_0)))$$
$$+ \sum_{k=1}^{K} n^{-k/2} E_{\theta_0} F(\Phi(Z + \sqrt{n} I^{1/2}(\theta_0)(\theta - \theta_0))) P_k(Z)$$
$$+ o(n^{-K/2}) h(n) + E_{\theta_0} R_n,$$

where the first expectation on the right-hand side is with respect to $Z \sim N(\mathbf{0}, I_d)$, and

$$h(n) = \int \frac{F(\Phi(z + \sqrt{n}I^{1/2}(\theta_0)(\theta - \theta_0)))}{1 + \|z\|^K} dz.$$

REMARK 1. In settings where our theorems for special cases do not apply, we can often obtain results by use of (2.10). This will be seen in Section 4. Moreover, it is seen that $h$ is integrable when $F(\Phi(Z + \sqrt{n}I^{1/2}(\theta_0)(\theta - \theta_0)))$ is.

PROOF OF PROPOSITION 2.1. Assumption JE gives that $W_o(\cdot|X^n)$ is approximated by $\Phi_{\mathbf{0}, I_d}(\cdot)$, or $W(\cdot|X^n)$ is approximated by $\Phi_{\hat{\theta}_n, (nI(\theta_0))^{-1}}(\cdot)$. Thus, the functional can be written as

$$F(W(\theta|X^n)) = F(\Phi_{\hat{\theta}_n, (nI(\theta_0))^{-1}}(\theta)) + R_n.$$

Taking expectations in $\theta_0$ and using Assumption EE gives

$$E_{\theta_0} F(W(\theta|X^n))$$

$$= \int F(\Phi_{u, (nI(\theta_0))^{-1}}(\theta))$$

$$\times \left( \phi_{\theta_0, (nI(\theta_0))^{-1}}(u) \right.$$

$$+ \sum_{k=1}^{K} n^{-k/2} P_k(\sqrt{n}I^{1/2}(\theta_0)(u - \theta_0)) \phi_{\theta_0, (nI(\theta_0))^{-1}}(u)$$

$$\left. + o(n^{-K/2}) \frac{|nI(\theta_0)|^{1/2}}{1 + \|\sqrt{n}I^{1/2}(\theta_0)(u - \theta_0)\|^K} \right) du + E_{\theta_0} R_n$$

$$= \int F(\Phi(z + \sqrt{n}I^{1/2}(\theta_0)(\theta - \theta_0))) \phi_d(z) dz$$

$$+ \sum_{k=1}^{K} n^{-k/2} \int F(\Phi(z + \sqrt{n}I^{1/2}(\theta_0)(\theta - \theta_0))) P_k(z) \phi_d(z) dz$$

$$+ o(n^{-(K/2)}) \int \frac{F(\Phi(z + \sqrt{n}I^{1/2}(\theta_0)(\theta - \theta_0)))}{1 + \|z\|^K} dz + E_{\theta_0} R_n. \quad \square$$

In examples we will see that $o(n^{-K/2})h(n)$ is often of lower order than $EF(\Phi(Z + \sqrt{n}I^{1/2}(\theta_0)(\theta - \theta_0)))$. Also, we observe the heuristic approximation

$$E[F(\Phi(Z + \sqrt{n}I^{1/2}(\theta_0)(\theta - \theta_0))) P_k(Z)]$$

$$\sim E[F(\Phi(Z + \sqrt{n}I^{1/2}(\theta_0)(\theta - \theta_0)))]E[P_k(Z)]$$
$$= E[F(\Phi(Z + \sqrt{n}I^{1/2}(\theta_0)(\theta - \theta_0)))]P_k(\sigma),$$

where $Z$ is a $N(\mathbf{0}, I_d)$ random vector, and $P_k(\sigma)$ is the expectation of $P_k(z)$ with powers $z^l$ replaced by $\sigma_l$, the $l$th moment of $N(\mathbf{0}, I_d)$. Taken together, these heuristics suggest that in many cases (2.10) gives

$$E_{\theta_0}F(W(\theta|X^n)) = EF(\Phi(Z + \sqrt{n}I^{1/2}(\theta_0)(\theta - \theta_0))) + \sum_{k=1}^{K} O(n^{-k/2}) + o(1).$$

**3. Asymptotics for expected values of functionals.** Proposition 2.1 was of general applicability. However, there are commonly occurring functionals that are worth examining in detail. When they depend on Johnson expandable quantities such as those in Theorem 2.1, we have a $K$-term expansion in powers of $n^{-j/2}$ on the "good" sets $S_n$. However, the coefficients depend on $X^n$. This is a problem because we want to take the expectation over the sample space for a functional of the posterior distribution. To get a closed form for these expectations, we must replace the empirical quantities in the coefficients in the expansion by their theoretical ones. Unfortunately, as noted in the remark after Theorem 2.1, such approximations are only accurate to order $o_p(1)$ unless more stringent hypotheses are proposed. Such hypotheses are hard to determine in part because the forms of the coefficients are generally unknown. Moreover, a posterior quantity must depend on the data, so replacing all the estimates with population values, if it could be done, defeats the purpose of using them. This is especially problematic when our goal is to obtain sample sizes. A final caveat is that we have tacitly been assuming that the expectation over the "bad" set $S_n^c$ will typically be small compared to that over the "good" set $S_n$, as noted in the Remark after Theorem 2.1, but we do not have a general closed form expression for it.

Taken together, these considerations mean we will only get a two-term expansion for the expectation, plus a remainder term

$$R_n' = E_{\theta_0}(F(W(\cdot|X^n)I_{S_n^c}),$$

which we have argued is asymptotically small enough, relative to the main approximation, that we can neglect it.

Theorems 3.1 and 3.2 below are extensions of results in [11], in which we have left the dimension of the parameter $d = 1$; cases with $d \geq 2$ are similar. Theorem 3.3 is more novel.

Let $\bar{\theta}$ be the posterior mean which often has the form $\bar{\theta} = s((1/n)\sum h(X_i)) + o_p(1/n)$. We use this in the first theorem because it is the right centering for posterior moments and is very close to the MLE. Note that in general we need to specify an estimator for planning purposes and that consistency of the MLE generally ensures that Bayes estimators are consistent; see [22]. Our first result is the following.

THEOREM 3.1. *Make all the assumptions in Section* 2, *in particular, those for Theorem* 2.1. *Also, assume Assumption* EE *for* $\bar{\theta}$ *in place of* $\hat{\theta}$. *Suppose* $\int |\theta|^r w(\theta)\, d\theta < \infty$ *and choose* $K, J \geq r$. *Then,*

$$(3.1) \quad E_{\theta_0} E_{W_o(\cdot|X^n)}[(\theta - \bar{\theta}_n)^r] = I^{-r/2}(\theta_0)\lambda_{rr} n^{-r/2} + o(n^{-r/2}) + R'_n,$$

*where* $\lambda_{rr} = 2^{r/2}\Gamma((r+1)/2)/\Gamma(1/2)$.

REMARK. In this case, the concern about using an approximation like $\hat{I}^{i/2}(\hat{\theta}_n) = I^{i/2}(\theta_0)(1 + o(1))$ for $i = 1, \ldots, r$ is built into Theorem 2.1: The scaling in the posterior by $I(\theta_0)$ and the laws of large number that are invoked to get $P_{\theta_0}(S_n) \to 0$ are enough for the expansions of posterior moments and percentiles.

PROOF OF THEOREM 3.1. Let $V_n = \sqrt{n} I^{1/2}(\theta_0)(\bar{\theta}_n - \theta_0)$. By Assumption EE for $V_n$, its density is

$$g_n(v) = \phi_d(v) + \sum_{k=1}^{K} n^{-k/2} P_k(v)\phi_d(v) + o(n^{-K/2})\frac{1}{1 + \|v\|^{K+2}}.$$

So we have

$$EV_n^r = \int v^r g_n(v)\, dv = \sigma_r + \sum_{k=1}^{K} n^{-k/2} P_{r,k}(\sigma) + o(n^{-K/2}),$$

where $\sigma$ is the vector of central moments from a $N(\mathbf{0}, I_d)$ as in (2.6) and the $o(n^{-K/2})$ comes from $o(n^{-K/2})\int v^r/(1 + \|v\|^{(K+2)})\, dv$. The integration is finite since $K \geq r$.

By using Assumption EE for both $\hat{\theta}_n$ and $\bar{\theta}_n$, we have

$$E_{\theta_0}(\hat{\theta}_n - \bar{\theta}_n) = I^{-1/2}(\theta_0)n^{-1/2} E_{\theta_0}(\alpha_n - V_n)$$

$$= I^{-1/2}(\theta_0)n^{-1/2} \sum_{k=1}^{K}(P_{1,k}(\sigma) - \bar{P}_{1,k}(\sigma))n^{-k/2} + o(n^{-K/2})$$

$$= O(n^{-1}),$$

where $\alpha_n = \sqrt{n} I^{1/2}(\hat{\theta} - \theta_0)$, the $P_{1,k}(\sigma)$'s are defined after (2.6), and the $\bar{P}_{1,k}(\sigma)$'s are their counterparts in the expansion for $f_{V_n}(\cdot)$. In general, for $m = 1, \ldots, r$, we have

$$(3.2) \qquad E_{\theta_0}(\hat{\theta}_n - \bar{\theta}_n)^m = O(n^{-(m+1)/2}).$$

Note $E_{\theta_0} E_{W_o(\cdot|X^n)}(\theta - \bar{\theta}_n)^r = E_{\theta_0} E_{W_o(\cdot|X^n)}(I_{S_n}(\theta - \bar{\theta}_n)^r) + R'_n$, and we only need to deal with the first of these terms. We omit the indicator $I_{S_n}$ for simplicity.

Assumption JE is satisfied by use of expression (2.2) in Theorem 2.1. Thus, for $i = 1, \ldots, r$ we have

$$(3.3) \quad E_{W_o(\cdot|X^n)}(I^{i/2}(\theta_0)(\theta - \hat{\theta}_n)^i) = \sum_{j=i}^{J} \lambda_{ij}(X^n) n^{-j/2} + O(n^{-(J+1)/2}),$$

on $\bigcap_{i=1}^{r} S_n(i)$ for $N \geq \max_{i=1}^{r} N_i$, where the $O(\cdot)$ is independent of $X^n$.

Now we can deal with the expectations $E_{\theta_0} E_{W_o(\cdot|X^n)} I^{i/2}(\theta_0)(\theta - \bar{\theta}_n)^i$, for $i = 1, \ldots, r$. Let $C(r, i)$ be the combination number of subsets of size $i$ from a set of size $r$. By (3.2) and (3.3), we have

$$E_{\theta_0} E_{W_o(\cdot|X^n)}(\theta - \bar{\theta}_n)^r$$
$$= E_{\theta_0} E_{W_o(\cdot|X^n)}((\theta - \hat{\theta}_n) + (\hat{\theta}_n - \bar{\theta}_n))^r$$
$$= (I(\theta_0))^{-r/2} E_{\theta_0} E_{W_o(\cdot|X^n)}(I^{r/2}(\hat{\theta}_n)(\theta - \hat{\theta}_n)^r)$$
$$\quad + \sum_{i=1}^{r} C(r,i) I^{-r/2}(\theta_0)$$
$$\qquad \times E_{\theta_0}[I^{(r-i)/2}(\theta_0)(\hat{\theta}_n - \bar{\theta}_n)^{r-i} E_{W_o(\cdot|X^n)}(I^{i/2}(\theta_0)(\theta - \hat{\theta}_n)^i)]$$
$$= I^{-r/2}(\theta_0) \lambda_{rr} n^{-r/2} + O(n^{-(r+1)/2})$$
$$\quad + \sum_{i=1}^{r} C(r,i) I^{-r/2}(\theta_0) O(n^{-(r-i+1)/2})$$
$$\qquad \times \left( \sum_{j=i}^{J} \lambda_{ij}(\theta_0) n^{-j/2} + O(n^{-(J+1)/2}) \right)$$
$$= I^{-r/2}(\theta_0) \lambda_{rr} n^{-r/2} + o(n^{-r/2}). \qquad \square$$

Now that we have an asymptotic form for functionals based on posterior moments, we turn to percentiles. Our result is the following.

THEOREM 3.2. *Make all the assumptions of Theorem 2.1 for some $J \geq 1$, and assume Assumption EE for some $K \geq 1$. Let $W^{-1}(\alpha|X^n)$ be the $\alpha$th quantile of $W(\cdot|X^n)$. Then we have*

$$(3.4) \quad E_{\theta_0} W^{-1}(\alpha|X^n) = \theta_0 + n^{-1/2} I^{-1/2}(\theta_0) \Phi^{-1}(\alpha) + o(n^{-1/2}) + R'(n).$$

PROOF. Let $\xi_\alpha$ be the $\alpha$th quantile of $\psi = \sqrt{n} I^{1/2}(\theta_0)(\theta - \hat{\theta}_n)$. That is,
$$\alpha = W_o(\psi \leq \xi_\alpha | X^n) = W(\theta \leq n^{-1/2} I^{-1/2}(\theta_0) \xi_\alpha + \hat{\theta}_n | X^n).$$

So, we get

$$W^{-1}(\alpha|X^n) = n^{-1/2}I^{-1/2}(\theta_0)\xi_\alpha + \hat{\theta}_n$$
(3.5)
$$= n^{-1/2}I^{-1/2}(\theta_0)\xi_\alpha + \theta_0 + n^{-1/2}I^{-1/2}(\theta_0)U_n,$$

where $U_n = \sqrt{n}I^{1/2}(\theta_0)(\hat{\theta}_n - \theta_0)$.

There is a function $\xi = \xi(\eta)$ which for any $\eta$ is a solution to $\Phi(\eta) = W_o(\xi(\eta)|X^n)$. So, given $\xi_\alpha$, we can backform to an $\eta_\alpha$ by defining the function $\xi(\cdot)$ to satisfy $\xi(\eta_\alpha) = \xi_\alpha$. Using this in (2.4) from Theorem 2.1, we get that $\xi_n(\alpha)$ satisfies Assumption JE, which we write as

$$\xi(\eta_\alpha) = \eta_\alpha + \sum_{j=1}^{J+1}\tau_j(\eta_\alpha)n^{-j/2}, \qquad n > N, X^n \in S_n,$$

where $\tau_{J+1}(\alpha)$ is the $O(n^{-(J+1)/2})$ remainder term, which is bounded in absolute value (a.s.). Using this in (3.5), we get

(3.6)
$$W^{-1}(\alpha|X^n) = n^{-1/2}I^{-1/2}(\theta_0)\left(\eta_\alpha + \sum_{j=1}^{J+1}\tau_j(\eta_\alpha)n^{-j/2}\right)$$
$$+ \theta_0 + n^{-1/2}I^{-1/2}(\theta_0)U_n,$$

for $n > N$ and $X^n \in S_n$.

By Assumption EE, we have

(3.7) $$E_{\theta_0}U_n = \sigma_1 + \sum_{k=1}^{K}n^{-k/2}P_{1,k}(\sigma) + o(n^{-K/2}),$$

in which we see $\sigma_1$ is the first moment of $N(0,1)$ and so is 0. Also, we have

$$\Phi(\eta_\alpha) = W_o(\xi(\eta_\alpha)|X^n) = W_o(\xi_\alpha|X^n) = \alpha,$$

so $\Phi^{-1}(\alpha) = \eta_\alpha$.

Finally, since $E_{\theta_0}W_n^{-1}(\alpha|X^n) = E_{\theta_0}I_{S_n}W_n^{-1}(\alpha|X^n) + R'_n$, we can take the expectation in (3.6), use (3.7), note that the $\tau_j(\eta_\alpha)$'s are bounded in $X^n$, collect terms and substitute for $\eta_\alpha$ to obtain

$$E_{\theta_0}W_n^{-1}(\alpha|X^n) = \theta_0 + n^{-1/2}I^{-1/2}(\theta_0)\Phi^{-1}(\alpha) + o(n^{-1/2}) + R'_n. \qquad \square$$

Next, we turn to derivatives of the posterior and more general posterior expectations. Denote the $r = (r_1,\ldots,r_d)$th derivative of $W(\theta|X^n)$ at $\theta$ by

$$W^{(r)}(\theta|X^n) = \frac{\partial^{|r|}}{\prod_{i=1}^{d}\partial\theta_i^{r_i}}W(\theta|X^n).$$

To express the first terms in the expansion for the expectation of a derivative of a posterior distribution, we need to define two sets of polynomials that arise when we differentiate expressions involving the normal density. The first is the set of Hermite polynomials: For a vector $i$ of length $d$, let $H_i(\cdot)$ be the $i$th Hermite polynomial defined by $H_{\mathbf{0}}(v) \equiv 1$ when $i = \mathbf{0}$ and by

$$D^{(i)}\phi(I^{1/2}(\theta_0)v) = H_i(v)\phi(I^{1/2}(\theta_0)v),$$

when $i \neq \mathbf{0}$. The second set of polynomials is particular to the use of Assumption JE for the posterior distribution. We define $\eta_j^{(r)}(\cdot)$ to be the polynomial given by

$$D^{(r)}[\phi(I^{1/2}(\theta_0)v)\gamma_j(I^{1/2}(\theta_0)v)] = \eta_j^{(r)}(v)\phi(I^{1/2}(\theta_0)v).$$

When we need to take expectations in the standard normal of products of polynomials $P(u)$ and $Q(u)$, we denote the polynomial of the normal moments by $P \circ Q$. That is, $EP(u)Q(u) \neq P(\sigma)Q(\sigma)$, but $EP(u)Q(u)$ is a polynomial in $\sigma$ which we denote $P \circ Q$. In this notation, we have the following.

THEOREM 3.3. *Assume Assumptions* JE *and* EE *for some* $J = K \geq 1$, *and that* $W(\theta|X^n)$ *has* $r = (r_1, \ldots, r_d)$th *derivative at* $\theta_0$ *with* $\min_i r_i \geq 1$. *Then*

$$
\begin{aligned}
E_{\theta_0} W^{(r)}(\theta_0|X^n) &= \frac{n^{|r|/2}|I(\theta_0)|^{1/2}}{(4\pi)^{d/2}} H_{r-1}\left(\frac{\sigma}{\sqrt{2}}\right) \\
&\quad + A_1 n^{(|r|-1)/2} + o(n^{(|r|-1)/2}) + R'_n,
\end{aligned}
\tag{3.8}
$$

*where* $r - 1 = (r_1 - 1, \ldots, r_d - 1)$, $H_{r-1}(\frac{\sigma}{\sqrt{2}})$ *is the expectation of* $H_{r-1}(v)$ *with powers* $v^s = v_1^{s_1} \cdots v_d^{s_d}$ *replaced by* $\sigma_s/(\sqrt{2})^{|s|}$ *and*

$$A_1 = \frac{1}{(4\pi)^{d/2}}\left(|I(\theta_0)|^{1/2} H_{r-1} \circ P_1\left(\frac{\sigma}{\sqrt{2}}\right) + \eta_1^{(r)}\left(\frac{\sigma}{\sqrt{2}}\right)\right).$$

PROOF. See the Appendix. □

If we set $r = (1, \ldots, 1)$ in Theorem 3.3, we get the posterior density. In fact, we can get the result for any partial derivative without the restriction $\min\{r_1, \ldots, r_d\} \geq 1$, by a similar technique. However, the computation of the coefficients becomes more involved. Also, in the Appendix we develop an asymptotic expansion for

$$E_{\theta_0}\left(\int h(\theta)w(\theta|X^n)\,d\theta\right),$$

where $h$ is a specified differentiable function; see (A.13). Such expansions may be helpful in sample size criteria derived from hypothesis testing optimality.

**4. Special cases.** Here, we examine four functionals encapsulating different sample size criteria taken from [23]. It will be seen that Proposition 2.1 and the results from Section 3 can be used to obtain closed form expressions for Bayesian sample sizes. To avoid repetition, we assume all the required conditions on the models are satisfied and just derive the corresponding formulae.

EXAMPLE 1. For the criterion APVC in [23], set

$$F(W(\cdot|X^n)) = \mathrm{Var}(\Theta|X^n)$$
$$= \int \theta'\theta W(d\theta|X^n) - \left(\int \theta W(d\theta|X^n)\right)'\left(\int \theta W(d\theta|X^n)\right).$$

By Theorem 3.1,

$$E_{\theta_0}F(W(\cdot|X^n)) = I^{-1}(\theta_0)\lambda_{22}n^{-1} + o(n^{-1}) + R'_n,$$

in which $\lambda_{22} = 2\Gamma(1+1/2)/\Gamma(1/2) = 1$, since $\Gamma(1+1/2) = 1/2\Gamma(1/2)$. Typically, $R'_n$ will be of smaller order than the main term, so for $\theta \in A$ with $\inf_{\theta \in A}|I(\theta)| > 0$, and prespecified $\varepsilon > 0$, the smallest sample size to achieve

$$|E_{\theta_0}\mathrm{Var}(\Theta|X^n)| \leq \varepsilon$$

is approximately given by

(APVC) $$n \geq \frac{1}{\varepsilon \inf_{\theta \in A} I(\theta)}.$$

A direct approach to this result by evaluating the terms in Proposition 2.1 can be done but seems to be quite involved.

EXAMPLE 2. For the criterion ACC in [23], set $F(W(\cdot|X^n)) = \int_{D_n} W(d\theta|X^n)$, in which $D_n$ is the HPD interval with length $l$ under the posterior distribution $W(\theta|X^n)$ and suppose $\theta$ is unidimensional. Unfortunately, our results in Section 3 do not apply, because, like the quantile example in the Introduction, the functional $F$ would have to depend on more than just the posterior.

However, we can still evaluate the terms in Proposition 2.1. The first term on the right-hand side of (2.6) is

$$EF(\Phi(Z + \sqrt{n}I^{1/2}(\theta_0)(\cdot - \theta_0)))$$

(4.2)
$$= \frac{\sqrt{n}I^{1/2}(\theta_0)}{\sqrt{2\pi}} E \int_{D'_n} e^{-(1/2)(Z+\sqrt{n}I^{1/2}(\theta_0)(\theta-\theta_0))^2} \, d\theta$$
$$= \frac{\sqrt{n}I^{1/2}(\theta_0)}{\sqrt{2\pi}} E \int_{D'_n} e^{-(1/2)nI(\theta_0)(\theta-\theta_0+Z/\sqrt{nI(\theta_0)})^2} \, d\theta.$$

From this, we see that $D'_n$ is of the form
$$D'_n = [\theta_0 - n^{-1/2}I^{-1/2}(\theta_0)Z - l/2, \theta_0 - n^{-1/2}I^{-1/2}(\theta_0)Z + l/2],$$
which is the HPD interval for $\theta$ under $\Phi(Z + \sqrt{n}I^{1/2}(\theta_0)(\cdot - \theta_0))$ of length $l$. Let $\eta = \sqrt{nI(\theta_0)}(\theta - \theta_0 + z/\sqrt{nI(\theta_0)})$. Then $\eta \sim N(0,1)$ and $D'_n = [-\sqrt{nI(\theta_0)}l/2 \leq \eta \leq \sqrt{nI(\theta_0)}l/2]$, so the right-hand side of (4.2) is

$$\frac{\sqrt{n}I^{1/2}(\theta_0)}{2\pi} \int \int_{D'_n} e^{-nI(\theta_0)/2(\theta-\theta_0+z/\sqrt{nI(\theta_0)})^2} \, d\theta \, e^{-(1/2)z^2} \, dz$$

$$= \frac{\sqrt{n}I^{1/2}(\theta_0)}{2\pi} \int \int_{[-\sqrt{nI(\theta_0)}l/2 \leq \eta \leq \sqrt{nI(\theta_0)}l/2]} e^{-\eta^2/2} \, d\eta \, e^{-(1/2)z^2} \, dz$$

$$= (2\Phi(\sqrt{nI(\theta_0)}l/2) - 1)\frac{1}{\sqrt{2\pi}} \int e^{-(1/2)z^2} \, dz$$

$$= 2\Phi(\sqrt{nI(\theta_0)}l/2) - 1.$$

As $n \to \infty$ this term tends to 1.

For large $n$, $D_n$ is of the form $[\bar{\theta}_n \pm l/2]$, where $\bar{\theta}_n$ is the posterior mean, and $\bar{\theta}_n \to \theta_0$ in $P_{\theta_0}$ probability. Also, we see that $F(\Phi(Z+\sqrt{n}I^{1/2}(\theta_0)(\cdot-\theta_0))$ is in fact independent of $Z$. Now, we have that

$$W([\tilde{\theta}_n \pm l/2]|X^n) \to 1$$

and

$$F(\Phi(Z + \sqrt{n}I^{1/2}(\theta_0)(\cdot - \theta_0))) = \Phi_{0,(nI(\theta_0))^{-1}}([\pm l/2]) \to 1,$$

also in $P_{\theta_0}$ probability. So, by the dominated convergence theorem, we have

$$E_{\theta_0}R_n = E_{\theta_0}(W([\tilde{\theta}_n \pm l/2]|X^n) - \Phi_{Z,(nI(\theta_0))^{-1}}([Z \pm l/2])) \to 0.$$

In the decomposition from Proposition 2.1(i), we see that (4.2) is the leading term and the other terms tend to zero. So, for given $0 < \alpha < 1$, the minimal $n$ to achieve

$$E_{\theta_0}F(W(\cdot|X^n)) = E_{\theta_0}\int_{D_n} W(d\theta|X^n) \geq 1 - \alpha$$

is approximately given by

$$2\Phi(\sqrt{nI(\theta_0)}l/2) - 1 \geq 1 - \alpha.$$

Equivalently, for $\theta \in A$ with $\inf_{\theta \in A} I(\theta) > 0$, we have

(ACC) $$n \geq \frac{4}{l^2 \inf_{\theta \in A} I(\theta)}\left[\Phi^{-1}\left(1 - \frac{\alpha}{2}\right)\right]^2,$$

where $\Phi(\cdot)$ is the distribution function of $N(0,1)$ and $\Phi^{-1}(\cdot)$ its inverse.

EXAMPLE 3. For the criterion ALC in [23], take $F(W(\cdot|X^n)) = W^{-1}_{\theta|X^n}(1-\alpha/2) - W^{-1}_{\theta|X^n}(\alpha/2)$, that is, suppose we require that the symmetric posterior quantiles be less than $l$ apart.

By Theorem 3.2,

$$E_{\theta_0} F(W(\cdot|X^n)) = \frac{1}{\sqrt{nI(\theta_0)}}(\Phi^{-1}(1-\alpha/2) - \Phi^{-1}(\alpha/2)) + o(n^{-1/2}).$$

So, for $\theta \in A$ with $\inf_{\theta \in A} I(\theta) > 0$, and given length $l$, the minimal $n$ to achieve

$$E_{\theta_0}(W^{-1}_{\theta|X^n}(1-\alpha/2) - W^{-1}_{\theta|X^n}(\alpha/2)) \leq l$$

is approximately given by

(ALC) $$n \geq \frac{1}{l^2 \inf_{\theta \in A} I(\theta)(\Phi^{-1}(1-\alpha/2) - \Phi^{-1}(\alpha/2))^2}.$$

Again, for completeness, we evaluate the terms in Proposition 2.1 directly. Let $\Phi_{Z,(nI(\theta_0))^{-1}}(\cdot)$ be the distribution function of $\phi_{Z,(nI(\theta_0))^{-1}}(\cdot)$ for given $Z$ and suppose $\theta$ is unidimensional. It is straightforward to see that

$$\Phi^{-1}_{Z,(nI(\theta_0))^{-1}}(\alpha/2) = Z + \frac{1}{\sqrt{nI(\theta_0)}}\Phi^{-1}(\alpha/2).$$

So, the first term in (2.10) is

$$E_{\theta_0} F(\Phi(Z + \sqrt{n}I^{1/2}(\theta_0)(\theta - \theta_0)))$$
$$= E_{\theta_0} F(\Phi_{Z,(nI(\theta_0))^{-1}}(\cdot))$$
$$= E_{\theta_0}(\Phi^{-1}_{Z,(nI(\theta_0))^{-1}}(1-\alpha/2) - \Phi^{-1}_{Z,(nI(\theta_0))^{-1}}(\alpha/2))$$
$$= E_{\theta_0}\left(\left(Z + \frac{1}{\sqrt{nI(\theta_0)}}\Phi^{-1}(1-\alpha/2)\right) - \left(Z + \frac{1}{\sqrt{nI(\theta_0)}}\Phi^{-1}(\alpha/2)\right)\right)$$
$$= \frac{1}{\sqrt{nI(\theta_0)}}(\Phi^{-1}(1-\alpha/2) - \Phi^{-1}(\alpha/2)),$$

as obtained above from Theorem 3.2.

Next, we deal with the remainder term in (2.6). In fact, it is enough to use (1.1), the two-term version of (2.6) avoiding nontrivial expansions entirely. Since we have

$$W(\sqrt{n}I^{1/2}(\theta_0)(\theta - \hat{\theta}_n)|X^n) \xrightarrow{d} N(\mathbf{0}, I_d),$$

we must have

$$W^{-1}_{\sqrt{n}I^{1/2}(\theta_0)(\theta - \hat{\theta}_n)|X^n}(\alpha) = \Phi^{-1}(\alpha) + o_p(1),$$

BAYESIAN SAMPLE SIZE 23$\forall\, 0 < \alpha < 1$. Equivalently,

$$W^{-1}_{(\theta|X^n)}(\alpha) = \hat{\theta}_n + \frac{1}{\sqrt{nI(\theta_0)}} \Phi^{-1}(\alpha) + o_p(n^{-1/2}).$$

So, we obtain

$$W^{-1}_{(\theta|X^n)}(\alpha) - \Phi^{-1}_{Z,(nI(\theta_0))^{-1}}(\alpha) = \hat{\theta}_n - Z + o_p(n^{-1/2}).$$

Since $E(Z) = \theta_0$, we can use Assumption EE to get

$$\begin{aligned} E_{\theta_0}(\hat{\theta}_n) &= \theta_0 + n^{-1/2} I^{-1/2}(\theta_0) E_{\theta_0}(\sqrt{n} I^{1/2}(\hat{\theta}_n - \theta_0)) \\ &= \theta_0 + n^{-1/2} I^{-1/2}(\theta_0) \\ &\quad \times \left( \int v \phi_d(v)\, dv \right. \\ &\quad \left. + \sum_{k=1}^{K} n^{k/2} \int v P_k(v) \phi_d(v)\, dv + o(n^{K+2}) \int \frac{v}{1 + \|v\|^{K+2}}\, dv \right) \\ &= \theta_0 + O(n^{-1}). \end{aligned}$$

Hence, with mild abuse of notation,

$$\begin{aligned} E_{\theta_0} R_n &= E_{\theta_0}(F(W(\cdot|X^n)) - F(\Phi_{Z,(nI(\theta_0))^{-1}}(\cdot))) \\ &= E_{\theta_0}(o_p(n^{-1/2})) = o(n^{-1/2}). \end{aligned}$$

EXAMPLE 4. For the *effect size* problem in [23], take $F(W(\cdot|X^n)) = \int_{\theta_1}^{\infty} W(d\theta|X^n)$. Here $\theta_1 < \theta_0$ and $\theta_1 < \hat{\theta}_n$ for large $n$. Our theorems do not apply, so we use Proposition 2.1. This gives

$$\begin{aligned} & EF(\Phi(Z + \sqrt{n} I^{1/2}(\theta_0)(\cdot - \theta_0))) \\ &\quad = E\left( \frac{\sqrt{nI(\theta_0)}}{\sqrt{2\pi}} \int_{\theta_1}^{\infty} e^{-(1/2)(Z + \sqrt{nI(\theta_0)}(\theta - \theta_0))^2}\, d\theta \right) \\ (4.3) &\quad = \frac{\sqrt{nI(\theta_0)}}{(\sqrt{2\pi})^2} \int \int_{\theta_1}^{\infty} e^{-(nI(\theta_0)/4)(\theta - \theta_0)^2}\, d\theta\, e^{-(z + (\sqrt{n}/2) I^{1/2}(\theta_0)(\theta - \theta_0))^2}\, dz \\ &\quad = \frac{\sqrt{nI(\theta_0)}}{\sqrt{2}\sqrt{2\pi}} \int_{\theta_1}^{\infty} e^{-(1/2)(nI(\theta_0)/2)(\theta - \theta_0)^2}\, d\theta \\ &\quad = 1 - \Phi\left( \frac{\sqrt{nI(\theta_0)}}{\sqrt{2}} (\theta_1 - \theta_0) \right). \end{aligned}$$

We see that (4.3) goes to 1 as $n$ increases (since $\theta_1 < \theta_0$). We show that the other terms are $O(n^{-1/2})$, so that (4.3) is the leading term.

In fact, since

$$F(\Phi(Z + \sqrt{n}I^{1/2}(\theta_0)(\cdot - \theta_0)))) = \frac{\sqrt{nI(\theta_0)}}{\sqrt{2\pi}} \int_{\theta_1}^{\infty} e^{-(1/2)(Z+\sqrt{nI(\theta_0)}(\theta-\theta_0))^2} \, d\theta$$

$$= 1 - \Phi(Z + \sqrt{nI(\theta_0)}(\theta_1 - \theta_0)),$$

which is bounded, for $J \geq 1$ we have that

$$\sum_{j=1}^{J} n^{-j/2} E[F(\Phi(Z + \sqrt{n}I^{1/2}(\theta_0)(\cdot - \theta_0)))) P_j(Z)] + o(n^{-1/2}) h(n)$$

$$= \sum_{j=1}^{J} n^{-j/2} O(1) E P_j(Z) + o(n^{-1/2}) h(n) = O(n^{-1/2}),$$

since the $EP_j(Z)$s are finite and $h(n) = o(1)$ by a similar evaluation as in (4.3).

For the remainder term, as in the proofs of the theorems, we only consider the "good" sets, omitting indicators on them. We have

$$E_{\theta_0} R_n = E_{\theta_0} \int_{\theta_1}^{\infty} d(W(\theta|X^n) - \Phi_{\hat{\theta}_n,(nI(\theta_0))^{-1}}(\theta))$$

$$= E_{\theta_0} \int_{\theta_1}^{\infty} \left( \sum_{j=1}^{J} n^{-j/2} n^{d/2} |I^{1/2}(\theta_0)| \right.$$

(4.4)
$$\times \phi_d(\sqrt{n}I^{1/2}(\theta_0)(\theta - \hat{\theta}_n)) \tilde{\gamma}_j(\sqrt{n}I^{1/2}(\theta_0)(\theta - \hat{\theta}_n))$$

$$\left. + n^{-(J-d+1)/2} |I^{1/2}(\theta_0)| \gamma_{J+1}^{(1)}(\sqrt{n}I^{1/2}(\theta_0)(\theta - \hat{\theta}_n)) \right) d\theta$$

$$= E_{\theta_0} \int_{\sqrt{n}I^{1/2}(\theta_0)(\theta_1-\hat{\theta}_n)}^{\infty} \left( \sum_{j=1}^{J} n^{-j/2} \phi_d(v) \tilde{\gamma}_j(v) + n^{-(K+1)/2} \gamma_{J+1}^{(1)}(v) \right) dv.$$

Since each term in (4.4) is integrable, expression (4.4) is bounded in absolute value by

$$\int \left( \sum_{j=1}^{K} n^{-j/2} \phi_d(v) |\tilde{\gamma}_j|(v) + n^{-(K+1)/2} |\gamma_{K+1}^{(1)}|(v) \right) dv = O(n^{-1/2}),$$

where, for a polynomial $P(\cdot)$, $|P|(v)$ is $P(v)$ with the coefficients and powers replaced by their absolute values.

So, for (4.3), for $\theta \in A = [a,b]$ with $\inf_{\theta \in A} I(\theta) > 0$, and given $0 < \alpha < 1$, the minimal $n$ achieving

$$E_{\theta_0} \int_{\theta_1}^{\infty} W(d\theta|X^n) \geq 1 - \alpha$$

is approximated by

$$\Phi\left(\frac{\sqrt{nI(\theta_0)}}{\sqrt{2}}(\theta_1 - \theta_0)\right) \leq \alpha,$$

which gives

(ES) $$n \geq \frac{2(\Phi^{-1}(\alpha))^2}{\inf_{\theta \in A}(\theta_1 - \theta)^2 I(\theta)}.$$

**5. Comparisons with exact results and numerical evaluations.** In this section we present some closed form expressions for the sample size criteria we have evaluated asymptotically. Then we turn to some numerical work. Both types of material suggest our asymptotic approximations are reasonable.

5.1. *Exact results.* In the case of the normal density with a conjugate normal prior we can obtain exact expressions from direct calculation for all four criteria we studied in Section 4. It is seen that our asymptotic expressions match these up to the stated error terms. More generally, only the (APVC) criterion, arguably the most popular of the four we have examined, can be calculated explicitly. We present two more examples, the Poisson($\theta$) with a Gamma($a, b$) prior and the Binomial($\theta$) with a Uniform([0, 1]) prior. Again, it is seen that our asymptotic expressions match the direct calculation expressions up to the stated order of error.

To begin the normal case, we record that, for $X|\theta \sim N(\theta, \sigma_0^2)$ and $\theta \sim N(\mu_0, \tau_0^2)$, we get $I(\theta_0) = \sigma_0^{-2}$ and $W(\theta|X^n) = N(\theta_n, \sigma_n^2)$ with

$$\theta_n = \frac{\overline{X} + \sigma_0^2 \mu_0/(n\tau_0^2)}{1 + \sigma_0^2/(n\tau_0^2)} \quad \text{and} \quad \sigma_n^2 = \frac{\sigma_0^2 \tau_0^2}{n\tau_0^2 + \sigma_0^2}.$$

Next we go through the four criteria in turn.

For the (APVC), the exact quantity is

$$E_{\theta_0}(\text{Var}(\theta|X^n)) = \text{Var}(\theta|X^n) = \frac{\sigma_0^2}{n + \sigma_0^2/\tau_0^2} = \frac{\sigma_0^2}{n} - \frac{\sigma_0^4/\tau_0^2}{n(n + \sigma_0^2/\tau_0^2)}.$$

If we choose $r = 2$, we have $\lambda_{22} = 2\Gamma(1 + 1/2)/\Gamma(1/2) = 1$, so by Theorem 3.1, we get

$$E_{\theta_0}(\text{Var}(\theta|X^n)) = \frac{\sigma_0^2}{n} + o(n^{-1}) + R'_n,$$

which matches up to the stated error.

For the (ALC), let $Z_n \sim N(\theta_n, \sigma_n^2)$. Then

$$\alpha = P(Z_n \leq W^{-1}(\alpha|X^n)|X^n) = P\left(\frac{Z_n - \theta_n}{\sigma_n} \leq \frac{W^{-1}(\alpha|X^n) - \theta_n}{\sigma_n}\Big|X^n\right),$$

so
$$\sigma_n^{-1}(W^{-1}(\alpha|X^n) - \theta_n) = \Phi^{-1}(\alpha) \quad \text{or} \quad W^{-1}(\alpha|X^n) = \theta_n + \sigma_n \Phi^{-1}(\alpha).$$

Since $\overline{X} \sim N(\theta, \sigma_0^2/n)$, we have
$$E_{\theta_0} W^{-1}(\alpha|X^n) = \frac{\theta_0 + \sigma_0^2 \mu_0/(n\tau_0^2)}{1 + \sigma_0^2/(n\tau_0^2)} + \frac{\sigma_0}{\sqrt{n + \sigma_0^2/\tau_0^2}} \Phi^{-1}(\alpha)$$
$$= \theta_0 + \frac{\sigma_0}{\sqrt{n}} \Phi^{-1}(\alpha) + o(n^{-1/2}).$$

By Theorem 3.2, we have
$$E_{\theta_0} W^{-1}(\alpha|X^n) = \theta_0 + \frac{\sigma_0}{\sqrt{n}} \Phi^{-1}(\alpha) + o(n^{-1/2}) + R'_n,$$

matching up to the stated error.

For the (ACC), we have $D_n = [\theta_n - l/2, \theta_n + l/2]$, and
$$E_{\theta_0} W([\theta_n \pm l/2]|X^n) = E_{\theta_0} \int_{[\theta_n - l/2, \theta_n + l/2]} \frac{1}{\sqrt{2\pi\sigma_n^2}} e^{-(1/(2\sigma_n^2))(\theta - \theta_n)^2} d\theta$$
$$= \int_{[-l/2, l/2]} \frac{1}{\sqrt{2\pi\sigma_n^2}} e^{-(1/(2\sigma_n^2))\alpha^2} d\alpha$$
$$= 2\Phi(\sigma_n l/2) - 1 = 2\Phi\left(\frac{\sqrt{n(1+\sigma_0^2)/(n\tau_0^2)}}{\sigma_0} \frac{l}{2}\right) - 1$$
$$= 2\Phi\left(\frac{\sqrt{n}}{\sigma_0} \frac{l}{2}\right) - 1 + o(1).$$

From Example 3, we have
$$E_{\theta_0} \int_{D_n} W(d\theta|X^n) d = 2\Phi\left(\frac{\sqrt{n}}{\sigma_0} \frac{l}{2}\right) - 1 + o(1),$$

matching up to the stated error.

For the effect size criterion, let $\theta_1 < \theta_0$. So, for large $n$, $\theta_n - \theta_1 \geq \varepsilon > 0$ (a.s.), $\sigma_n^{-1} = O(n^{1/2})$ and $\theta_n - \overline{X} = O(n^{-1})$, giving $\sigma_n^{-1}(\theta_n - \overline{X}) = O(n^{-1/2})$. Using this in the functional gives
$$E_{\theta_0} \int_{\theta_1}^{\infty} W(d\theta|X^n)$$
$$= E_{\theta_0} \int_{\theta_1}^{\infty} \frac{1}{\sqrt{2\pi}\sigma_n} e^{-(\theta-\theta_n)^2/(2\sigma_n^2)} d\theta$$
$$= E_{\theta_0} \int_{\sigma_n^{-1}(\theta_1 - \theta_n)}^{\infty} \frac{1}{\sqrt{2\pi}} e^{-\alpha^2/2} d\alpha$$

$$= E_{\theta_0} \int_{\sigma_n^{-1}(\theta_1-\overline{X})}^{\infty} \frac{1}{\sqrt{2\pi}} e^{-\alpha^2/2}\, d\alpha + E_{\theta_0} \left| \int_{\sigma_n^{-1}(\theta_1-\overline{X})}^{\sigma_n^{-1}(\theta_1-\theta_n)} \frac{1}{\sqrt{2\pi}} e^{-\alpha^2/2}\, d\alpha \right|$$

$$= E_{\theta_0} \int_{\sigma_n^{-1}(\theta_1-\overline{X})}^{\infty} \frac{1}{\sqrt{2\pi}} e^{-\alpha^2/2}\, d\alpha + O(n^{-1/2}).$$

Since

$$\sigma_n^{-1} - \frac{\sqrt{n}}{\sigma_0} = \frac{1}{\tau_0^2(\sqrt{n\tau_0^2 + \sigma_0^2}/(\sigma_0\tau_0) + \sqrt{n}/\sigma_0)} = O(n^{-1/2}),$$

so $\sigma_n^{-1}(\theta_1 - \overline{X}) - \sqrt{n}/\sigma_0(\theta_1 - \overline{X}) = O(n^{-1/2})$ (a.s.), the last expression for the functional is

$$E_{\theta_0} \int_{\sqrt{n}/\sigma_0(\theta_1-\overline{X})}^{\infty} \frac{1}{\sqrt{2\pi}} e^{-\alpha^2/2}\, d\alpha$$

$$+ E_{\theta_0} \left| \int_{\sigma_n^{-1}(\theta_1-\overline{X})}^{\sqrt{n}/\sigma_0(\theta_1-\overline{X})} \frac{1}{\sqrt{2\pi}} e^{-\alpha^2/2}\, d\alpha \right| + O(n^{-1/2})$$

$$= E_{\theta_0} \int_{\sqrt{n}/\sigma_0(\theta_1-\overline{X})}^{\infty} \frac{1}{\sqrt{2\pi}} e^{-\alpha^2/2}\, d\alpha + O(n^{-1/2}) + O(n^{-1/2})$$

$$= \int_{-\infty}^{\infty} \int_{\theta_1}^{\infty} \frac{1}{\sqrt{2\pi}} \frac{\sqrt{n}}{\sigma_0} e^{-n(\theta-x)^2/(2\sigma_0^2)} \frac{1}{\sqrt{2\pi}} \frac{\sqrt{n}}{\sigma_0} e^{-n(x-\theta_0)^2/(2\sigma_0^2)}\, d\theta\, dx$$

$$+ O(n^{-1/2})$$

$$= \int_{\theta_1}^{\infty} \frac{1}{\sqrt{4\pi}} \frac{\sqrt{n}}{\sigma_0} e^{-n(\theta-\theta_0)^2/4\sigma_0^2} \int_{-\infty}^{\infty} \frac{1}{\sqrt{\pi}} \frac{\sqrt{n}}{\sigma_0} e^{-n(x-(\theta+\theta_0)/2)^2/\sigma_0^2}\, dx\, d\theta$$

$$+ O(n^{-1/2})$$

$$= \int_{\theta_1}^{\infty} \frac{1}{\sqrt{4\pi}} \frac{\sqrt{n}}{\sigma_0} e^{-n(\theta-\theta_0)^2/(4\sigma_0^2)}\, d\theta + O(n^{-1/2})$$

$$= \int_{(\sqrt{n}/(\sqrt{2}\sigma_0))(\theta_1-\theta_0)}^{\infty} \frac{1}{\sqrt{2\pi}} e^{-\alpha^2/2}\, d\alpha + O(n^{-1/2})$$

$$= 1 - \Phi\left(\sqrt{\frac{n}{2}} \frac{\theta_1 - \theta_0}{\sigma_0}\right) + O(n^{-1/2}).$$

From Example 4, we have that

$$E_{\theta_0} \int_{\theta_1}^{\infty} W(d\theta|X^n) = 1 - \Phi\left(\sqrt{\frac{n}{2}} \frac{\theta_1 - \theta_0}{\sigma_0}\right) + o(1),$$

again matching up to the stated error. In this case, the exact expression gave slightly stronger control of the error.

Next, we turn to two other examples for the (APVC). Of the four criteria, only the (APVC) is simple enough that it can be obtained in closed form in some cases.

Let $X|\theta \sim \text{Poisson}(\theta)$, and suppose $\theta \sim G(a,b)$, the Gamma distribution with $a,b$ known. Let $S_n = \sum_{i=1}^n X_i$. Then, by standard results, $\theta|X^n \sim G(a+n, b+S_n)$, with $E(\theta|X^n) = (b+S_n)/(a+n)$, $\text{Var}(\theta|X^n) = (b+S_n)/(a+n)^2$, $I(\theta_0) = 1/\theta_0$, and $E_{\theta_0}(S_n) = n\theta_0$.

So, the expected posterior variance is

$$E_{\theta_0}(\text{Var}(\theta|X^n)) = \frac{b+n\theta_0}{(a+n)^2} = \frac{\theta_0}{n} + \frac{b-\theta_0}{(n+a)^2} - \frac{a}{n(n+a)},$$

and by Theorem 3.1, the approximation is

$$E_{\theta_0}(\text{Var}(\theta|X^n)) = \frac{\theta_0}{n} + o(n^{-1}) + R'_n.$$

As in the normal case, the two match up to the stated error.

Now, let $X|\theta \sim \text{Binomial}(\theta)$ with $\theta \sim U(0,1)$. Setting $S_n = \sum_{i=1}^n X_i$, standard results give that $\theta|X^n \sim \text{Beta}(S_n + 1, n + 1 - S_n)$, with $E(\theta|X^n) = (S_n+1)/(n+2)$, $\text{Var}(\theta|X^n) = (nS_n - S_n^2 + n + 1)/[(n+2)^2(n+3)]$, $I(\theta_0) = 1/[\theta_0(1-\theta_0)]$, $E_{\theta_0}(S_n) = n\theta_0$ and $E_{\theta_0}(S_n^2) = n\theta_0(1-\theta_0) + n^2\theta_0^2$.

The expected posterior variance is

$$E_{\theta_0}(\text{Var}(\theta|X^n))$$
$$= \frac{n^2\theta_0 - n\theta_0 - n(n-1)\theta_0^2 + n + 1}{(n+2)^2(n+3)}$$
$$= \frac{\theta_0(1-\theta_0)}{n} - \frac{3\theta_0(1-\theta_0)}{n(n+3)} + \frac{1-\theta_0(1-\theta_0)}{(n+2)(n+3)} - \frac{2\theta_0(1-\theta_0)+1}{(n+2)^2(n+3)}.$$

By Theorem 3.1 our approximation is

$$E_{\theta_0}(\text{Var}(\theta|X^n)) = \frac{\theta_0(1-\theta_0)}{n} + o(n^{-1}) + R'_n.$$

As before, the two agree.

The agreement between the asymptotics and the closed form expressions suggests that in the other examples the discrepancy between the two will be small. Indeed, all of the criteria are derived from posteriors and posterior objects which can be approximated as well as desired by taking enough terms in the expansions. That is, optimizing the asymptotic expression obtained by using more terms will give any desired degree of accuracy. We suggest this will only be needed in extreme cases when the coefficients in the neglected higher-order terms are so large, possibly because of the range of the set in the parameter space, that they overwhelm the lower-order terms.

5.2. *Numerical evaluations.* Fundamentally, the class of quantities we want to evaluate is of the form $G = E_\theta F_\varepsilon(W(\cdot|X^n))$, where $F$ represents the inference objective and $\varepsilon$ summarizes how well it must be met. To begin, we

TABLE 1
*Exact vs. asymptotic: Normal–Normal*

| Parameter | $n$ | (**APVC**): $G$, $G^*$ | (**ALC**): $G$, $G^*$ | (**ACC**): $G$, $G^*$ |
|---|---|---|---|---|
| $\eta_1$ | 10  | 0.0187 (0.0200) | 0.2591 (0.2674) | 0.1449 (0.1403) |
|          | 30  | 0.0065 (0.0067) | 0.3617 (0.3657) | 0.2431 (0.2405) |
|          | 50  | 0.0039 (0.0040) | 0.3934 (0.3960) | 0.3093 (0.3074) |
|          | 100 | 0.0020 (0.0020) | 0.4250 (0.4264) | 0.4251 (0.4238) |
| $\eta_2$ | 10  | 0.2308 (0.2500) | 4.0944 (4.1776) | 0.3972 (0.3829) |
|          | 30  | 0.0811 (0.0833) | 4.4911 (4.5252) | 0.6200 (0.6135) |
|          | 50  | 0.0492 (0.0500) | 4.6106 (4.6322) | 0.74040 (0.7364) |
|          | 100 | 0.0248 (0.0250) | 4.7286 (4.7399) | 0.8877 (0.8862) |
| $\eta_3$ | 10  | 1.6071 (1.8000) | 22.3791 (22.7932) | 0.6759 (0.6485) |
|          | 30  | 0.5769 (0.6000) | 23.5583 (23.7259) | 0.9002 (0.8934) |
|          | 50  | 0.3516 (0.3600) | 23.9075 (24.0131) | 0.9650 (0.9628) |
|          | 100 | 0.1779 (0.1800) | 24.2470 (24.3022) | 0.9970 (0.9968) |

present computations for two simple cases in which $G$ can be obtained from the closed form expressions in Section 5.1. We compare selected values of $G$ with the corresponding approximations $G^*$ from our asymptotic formulae. We look at expected values of functionals, rather than fix $\varepsilon$'s and find optimal sample sizes, to make it easy to compare these first two simple cases with a more complicated third case.

Table 1 gives the exact $G$ and approximate $G^*$ (in brackets) numerical results for the normal likelihood and normal prior example given in Section 5.1. We have set $\eta = (\theta_0, \mu_0, \sigma_0^2, \tau_0^2)$ and chosen $\eta_1 = (0.5, 0.25, 0.20, 0.30)$, $\eta_2 = (5.0, 3.5, 2.5, 3.0)$ and $\eta_3 = (25, 20, 18, 15)$; the values of $n$ are as indicated. The confidence level for (ALC) is $\alpha = 0.05$; for (ACC), we set $l = \theta_0/10$. (We omitted results for the effect size problem because the exact and the approximate quantities have the same first-order term and the higher-order terms are hard to get explicitly.)

It is seen that as $n$ increases the values of the (APVC) functional decrease, while the values for (ALC) and (ACC) increase. This is expected from the interpretation of the functionals. For each choice of $\eta$ and criterion, it is seen that the error decreases as $n$ increases; that is, the difference between $\hat{G}$ and $G^*$ gets smaller as $n$ gets larger. It is important to note that, as the numerical value of $G$ changes, it is closely tracked by our approximation.

Less routine examples are the Poisson$(\theta)$ likelihood with a Gamma$(a, b)$ prior and a binomial $(\theta)$ likelihood with a Uniform$[0, 1]$ prior. For the Poisson–Gamma, we set $\eta = (\theta_0, a, b)$ and for the Binomial–Uniform we set $\eta = \theta_0$.

Table 2 shows the values for (APVC) from $G$ and $G^*$ for $\eta_1 = (0.5, 2.5, 3.5)$, $\eta_2 = (1.6, 8, 7.5)$ and $\eta_3 = (1.5, 10, 12)$. For the Binomial–Uniform, we set $\eta_1 = 0.20$, $\eta_2 = 0.5$ and $\eta_3 = 0.75$.

As in Table 1, both the error of approximation and the numerical values decrease as $n$ increases for both prior likelihood pairs. For the Poisson–Gamma case, it is seen that the values for $\eta_2$ and $\eta_3$ are close because their $\theta$'s are close. The prior has a smaller effect. For the Binomial–Uniform with constant prior, it is seen that the symmetry of the Binomial makes the values for $\eta_1$ and $\eta_2$ close.

Next, we turn to an example in which a closed form for $G$ does not exist. We will approximate $G$ by $\hat{G}$ obtained from simulations and compare this to $G^*$ again obtained from our asymptotic expressions. To clarify the comparison in Table 3, observe that, in a world of infinite resources, we would generate $m$ IID $X^n$'s from $p_\theta$, find $W(\cdot|X^n = x^n)$ for each of the $x_j^n$'s, evaluate $\hat{G}(\theta, \varepsilon, W, n, m) = (1/m) \sum_{j=1}^m F_\varepsilon(W(\cdot|X^n = x_j^n))$ and report $\hat{G} = \hat{G}(\theta, \varepsilon, W, n, m)$ as an approximation to $G = G(\theta, \varepsilon, W, n)$. Ideally, we would use a large enough $m$ that dependence on it could be neglected and $W$ would be replaced by the hyperparameters, say, $\alpha$, that specify it. That is, we will have

$$(5.1) \qquad \hat{G}(\theta, \varepsilon, \alpha, n, m) \approx G(\theta, \varepsilon, \alpha, n),$$

so we can obtain minimizing values of $n = n(\theta, \varepsilon, \alpha)$ from $\hat{G}$. In fact, we want a maximin solution

$$(5.2) \qquad n_{Mm}(\varepsilon) = \max_{\theta \in K, \alpha \in A} n(\theta, \varepsilon, \alpha),$$

in which $K$ and $A$ are compact sets. However, direct evaluation of $n_{Mm}(\varepsilon)$ is computationally demanding: It requires, for each specified $\varepsilon$, $\theta$ and $\alpha$, evaluating $E_\theta F_\varepsilon(W(\cdot|X^n))$ for many values of $n$ so one can select the smallest $n$ that satisfies the criterion.

As in the first two cases, rather than evaluating (5.2), we compute, for some choices of $n$, the empirical posterior functional $\hat{G}(\theta, \varepsilon, \alpha, n, m)$, which BAYESIAN SAMPLE SIZE 31can be regarded as a good enough approximation to $G(\theta, \varepsilon, \alpha, n)$ for large $m$. We also compute our asymptotic approximation, $G^*$. In effect, we have assumed (5.1) by choosing $m$ large enough and then compared $\hat{G}(\theta, \varepsilon, \alpha, n)$ to $G^*(\theta, \varepsilon, \alpha, n)$. Thus, Table 3 gives $G^*$ and $\hat{G}$ for several choices of $\theta$, $\varepsilon$, $\alpha$ and $n$, for various functionals $F$.

TABLE 2
*Exact vs. asymptotic: Non-Normal*

| $\eta$ \ $n$ | 10 | 30 | 50 | 100 |
|---|---|---|---|---|
| | | Poisson–Gamma | | |
| $\eta_1$ | 0.0544 (0.0500) | 0.0175 (0.0167) | 0.0103 (0.0100) | 0.0051 (0.0050) |
| $\eta_2$ | 0.0725 (0.1600) | 0.0384 (0.0533) | 0.0260 (0.0320) | 0.0144 (0.0160) |
| $\eta_3$ | 0.0675 (0.1500) | 0.0356 (0.0500) | 0.0242 (0.0300) | 0.0134 (0.0150) |
| | | Binomial–Uniform | | |
| $\eta_1$ | 0.0136 (0.0160) | 0.0050 (0.0053) | 0.0031 (0.0032) | 0.0016 (0.0016) |
| $\eta_2$ | 0.0179 (0.0250) | 0.0073 (0.0083) | 0.0046 (0.0050) | 0.0024 (0.0025) |
| $\eta_3$ | 0.0149 (0.0188) | 0.0057 (0.0063) | 0.0036 (0.0038) | 0.0018 (0.0019) |

Our argument is that the approximations $G^*$ are close to the corresponding $\hat{G}$'s for a variety of points $(\theta, \varepsilon, \alpha, n)$ and, therefore, it is reasonable to use sample sizes obtained from $G^*$ as approximations to the sample sizes one would get from optimizing $G$ directly. The values given for the $\hat{G}$ and $G^*$ given in the tables support this contention.

Thus, we evaluated a nonconjugate, nonclosed form example. In this case, the $G$ could not be found as in Section 5.1; we are forced to use $\hat{G}$. To provide a real test of the asymptotics, take the likelihood to be Exponential$(x|\theta) = \theta e^{-\theta x}$ with a Beta(3/2, 3/2) prior having density $\beta(\theta|3/2, 3/2) \propto \sqrt{\theta(1-\theta)}$ on [0, 1]. It is seen that this example is far from the normal prior, normal likelihood setting, so its relation to the asymptotic normality used to derive our expressions is not close.

Since $G$ is an expected value of a functional of the posterior, we generate $m = 1000$ IID data sets of size $n$ for several values of $n$, $X_1^n, \ldots, X_m^n$, from an Exponential$(x|\theta)$. For each $X_j^n$, $j = 1, \ldots, m$, we draw outcomes from $W(\cdot|X_j^n)$ by Markov chain Monte Carlo, compute $F(W(\cdot|X_j^n))$ from the empirical posterior distribution, and approximate $E_\theta F(W(\cdot|X^n))$ by $(1/m) \sum_{j=1}^m F(W(\cdot|X_j^n))$.

For several values of $\theta$ taken as true, $n$ as a potential sample size, and each of three criteria, we give the empirical value, $\hat{G}$, and its asymptotic approximation using our technique $G^*$ in brackets in Table 3. The expected

TABLE 3
*Empirical vs. asymptotic: Non-Normal*

| $\theta_0$ | $n$ | $E_{\theta_0}(\text{Var}(\theta|X^n))$ | $E_{\theta_0}(\text{HPD})$ | $E_{\theta_0}(ALC)$ |
|---|---|---|---|---|
| 0.25 | 10  | 0.0116 (0.0062) | [0.1475, 0.5388] ([0.1633, 0.4732]) | 0.3912 (0.3099) |
|      | 30  | 0.0031 (0.0021) | [0.1742, 0.3826] ([0.1803, 0.3592]) | 0.2084 (0.1789) |
|      | 50  | 0.0018 (0.0012) | [0.1884, 0.3483] ([0.1939, 0.3325]) | 0.1599 (0.1386) |
|      | 100 | 0.0008 (0.0006) | [0.2017, 0.3123] ([0.2055, 0.3035]) | 0.1106 (0.0980) |
| 0.50 | 10  | 0.0238 (0.0250) | [0.2703, 0.8399] ([0.2320, 0.8518]) | 0.5696 (0.6198) |
|      | 30  | 0.0107 (0.0083) | [0.3387, 0.7273] ([0.3409, 0.6988]) | 0.3886 (0.3578) |
|      | 50  | 0.0068 (0.0050) | [0.3727, 0.6832] ([0.3798, 0.6570]) | 0.3105 (0.2772) |
|      | 100 | 0.0034 (0.0025) | [0.4020, 0.6208] ([0.4084, 0.6044]) | 0.2188 (0.1960) |
| 0.75 | 10  | 0.0348 (0.0562) | [0.3738, 0.9467] ([0.2135, 1.1432]) | 0.5729 (0.9297) |
|      | 30  | 0.0140 (0.0187) | [0.4986, 0.9368] ([0.4556, 0.9923]) | 0.4382 (0.5368) |
|      | 50  | 0.0102 (0.0112) | [0.5511, 0.9282] ([0.5349, 0.9506]) | 0.3771 (0.4158) |
|      | 100 | 0.0059 (0.0056) | [0.5988, 0.8988] ([0.5997, 0.8937]) | 0.3000 (0.2940) |

HPD is a proxy for (ACC): In the average coverage criterion, we fix $\ell$ and find the $n$ making the coverage probability of the HPD set of length less than $\ell$ at least $1 - \alpha$. Here, the $E(HPD)$ represents the $\ell$ for coverage 0.95 for the approximate HPD interval centered at the posterior mean.

It is seen that the expected (APVC) and (ALC) decreases as $n$ increases, as does the error of approximation. Likewise, the expected HPD length decreases, as does the error of approximation as $n$ increases. When $n = 10$, the approximation can be poor with errors often over 25% of the true value; this may be due to the $m$ or $n$ being too small or due to convergence problems in the Markov chain Monte Carlo. At the other end, $n = 100$ gives good approximation in absolute and relative senses, suggesting the size of $m$ is not the problem. Overall, in highly nonnormal and nonconjugate settings, our approximation may not give satisfactory results unless $n$ is moderate, say, over 30.

We comment that the effect size criterion involves the mean posterior quantiles, so we expect our formulae to give results similar to those for $E_{\theta_0}(HPD)$, for which reason we omitted its presentation here.

**6. Final remarks.** Overall, we have argued that simple, asymptotically valid inequalities can be derived so that Bayesian sample sizes can be readily determined essentially as easily as in the frequentist case. We have done this for four sample size criteria taken from the established literature.

Apart from this contribution, we have several observations.

First, integrating our approximations for (1.1) over $\theta_0$ gives expressions for use in pre-posterior Bayesian calculations where the expectations are taken with respect to the mixture density. That is, because $F(W(\cdot|X^n))$ does not depend on the parameter explicitly, the expectation with respect to the mixture is $E_m F(W(\cdot|X^n)) = \int_\Theta E_\theta F(W(\cdot|X^n)) w(\theta) \, d\theta$, and our asymptotic expressions will apply to the argument of the integral. Our results are slightly stronger than necessary for evaluating marginal probabilities.

Second, although we have not done it here, we suggest that, as ever, sensitivity analyses should be used to ensure the sample sizes obtained from any one method are robust against deviations of the prior, likelihood and loss function (if one exists) from the nominal choices used to obtain the sample sizes. Robustness against similar choices of sample size criterion would also be desirable.

Finally, we anticipate that examining functionals of posteriors may be a step toward unifying the three cases described in the Introduction. Decision theoretic procedures implicitly rest on the posterior risk which can be regarded as a functional of the posterior. Evidentiary procedures usually devolve to posterior probabilities which can likewise be regarded as functionals of the posterior—we suggest formulae for these at the end of the Appendix. And, finally, purely Bayes criteria that focus on credibility sets also express properties of credibility sets in terms of the posterior. It may be that a suitably general treatment of functionals of the posterior will include all these as special cases of one unified formalism.

## APPENDIX

Here, we prove Theorem 3.3 and compare it with the expansions for two functionals in [5]. As a final point, we note how to use our techniques to get an asymptotic expansion for a functional that is the expectation of a posterior mean of a function of the parameter.

PROOF OF THEOREM 3.3. We need to approximate $E_{\theta_0}(I_{S_n} W^{(r)}(\theta_0|X^n))$; for simplicity of notation, we omit the $I_{S_n}$.

First, for $1 \leq j \leq J$, the $\gamma_j(\sqrt{n} I^{1/2}(\theta_0)(\theta - \hat{\theta}_n), X^n)$'s are polynomials and, hence, differentiable. As in Assumption JE, the remainder term is

$$\gamma_{J+1}(\sqrt{n} I^{1/2}(\theta_0)(\theta_0 - \hat{\theta}_n), X^n) n^{-(J+1)/2}$$
$$= W(\theta_0|X^n) - \Phi_{\hat{\theta}_n, I^{-1}(\theta_0)/n}(\theta_0)$$
$$- \sum_{j=1}^{J} n^{-j/2} \phi_d(\sqrt{n} I^{1/2}(\theta_0)(\theta_0 - \hat{\theta}_n)) \gamma_j(\sqrt{n} I^{1/2}(\theta_0)(\theta_0 - \hat{\theta}_n), X^n),$$

$$n > N, X^n \in S_n.$$

So $\gamma_{K+1}(\cdot, X^n)$ has $r$th derivative whenever $W(\cdot|X^n)$ does.

To control the expectation of $W^{(r)}(\theta|X^n)$, we replace the $\gamma_j(\cdot, X^n)$'s by the $\gamma_j(\cdot)$'s. That is, by the boundedness of the $\gamma_j(\cdot, X^n)$'s, and the a.s. convergence of $\hat{\theta}_n$ and the $I_r(\hat{\theta}_n)$'s to $\theta_0$ and the $I_r(\theta_0)$'s, we have

$$\gamma_j(\cdot, X^n) = \gamma_j(\cdot)(1 + o_p(1)),$$

for $j = 1, \ldots, J+1$, where the $o_p(1)$ may depend on $j$, but is independent of $\theta$. So we have

$$W(\theta_0|X^n) = \Phi_{\hat{\theta}_n, I^{-1}(\theta_0)/n}(\theta_0)$$
$$+ \sum_{j=1}^{J} n^{-j/2} \phi_d(\sqrt{n} I^{1/2}(\theta_0)(\theta_0 - \hat{\theta}_n))$$
(A.1)
$$\times \gamma_j(\sqrt{n} I^{1/2}(\theta_0)(\theta_0 - \hat{\theta}_n))(1 + o_p(1))$$
$$+ n^{-(J+1)/2} \gamma_{J+1}(\sqrt{n} I^{1/2}(\theta_0)(\theta_0 - \hat{\theta}_n))(1 + o_p(1)).$$

Next, we convert (A.1) into a form to which Assumption EE can be applied. We begin to deal with derivatives of the first term by noting

$$\left.\frac{\partial^{|r|} \Phi_{\hat{\theta}_n, I^{-1}(\theta_0)/n}(\theta)}{\partial \theta^r}\right|_{\theta=\theta_0} = \left.\frac{\partial^{|r-1|} \phi_{\hat{\theta}_n, I^{-1}(\theta_0)/n}(\theta)}{\partial \theta^{r-1}}\right|_{\theta=\theta_0}.$$

Next, let $I_i^{1/2}(\theta_0)$ be the $i$th column of $I^{1/2}(\theta_0)$, and $1_i = (0, \ldots, 0, 1, 0, \ldots, 0)$ be the $d$-vector with the $i$th component 1 and all other components zero. For the first derivative with respect to the $i$th component of $\theta$ we have

$$\frac{\partial \phi_{\hat{\theta}_n, I^{-1}(\theta_0)/n}(\theta)}{\partial \theta_i}$$

$$= |nI(\theta_0)|^{1/2} \frac{\partial \phi(\sqrt{n} I^{1/2}(\theta_0)(\theta - \hat{\theta}_n))}{\partial \theta_i}$$

$$= \left( \frac{|nI(\theta_0)|^{1/2} \partial \phi(\sqrt{n} I^{1/2}(\theta_0)(\theta - \hat{\theta}_n))}{\partial [\sqrt{n} I^{1/2}(\theta_0)(\theta - \hat{\theta}_n)]} \right)' \frac{\partial}{\partial \theta_i} (\sqrt{n} I^{1/2}(\theta_0)(\theta - \hat{\theta}_n))$$

$$= n^{(d+1)/2} |I(\theta_0)|^{1/2} I_i^{1/2}(\theta_0)(\sqrt{n} I^{1/2}(\theta_0)(\theta - \hat{\theta}_n)) \phi(\sqrt{n} I^{1/2}(\theta_0)(\theta - \hat{\theta}_n))$$

$$= n^{(d+1)/2} |I(\theta_0)|^{1/2} H_{1_i}(\sqrt{n} I^{1/2}(\theta_0)(\theta - \hat{\theta}_n)) \phi(\sqrt{n} I^{1/2}(\theta_0)(\theta - \hat{\theta}_n)).$$

So, by an induction argument we have

$$\begin{aligned}
\left. \frac{\partial^{|r-1|} \phi_{\hat{\theta}_n, I^{-1}(\theta_0)/n}(\theta)}{\partial \theta^{r-1}} \right|_{\theta=\theta_0} \\
= n^{|r|/2} |I(\theta_0)|^{1/2} H_{r-1}(\sqrt{n} I^{1/2}(\theta_0)(\theta_0 - \hat{\theta}_n)) \\
\times \phi(\sqrt{n} I^{1/2}(\theta_0)(\theta_0 - \hat{\theta}_n)),
\end{aligned} \quad (A.2)$$

in which we have simplified by using $(d + |r - 1|)/2 = |r|/2$.

Using (A.2) in the first term, and recalling the definition of the $\eta_j^{(r)}$ in the second term, the $r$th derivative of (A.1) becomes

$$W^{(r)}(\theta_0 | X^n)$$

$$(A.3) \quad = n^{|r|/2} |I(\theta_0)|^{1/2} H_{r-1}(\sqrt{n} I^{1/2}(\theta_0)(\theta_0 - \hat{\theta}_n))$$

$$\times \phi(\sqrt{n} I^{1/2}(\theta_0)(\theta - \hat{\theta}_n))$$

$$+ \sum_{j=1}^{J} n^{-j/2} n^{|r|/2} \eta_j^{(r)}(\sqrt{n} I^{1/2}(\theta_0)(\theta_0 - \hat{\theta}_n))$$

$$(A.4) \quad \times (1 + o_p(1)) \phi(\sqrt{n} I^{1/2}(\theta_0)(\theta_0 - \hat{\theta}_n))$$

$$(A.5) \quad + n^{-(J+1)/2} n^{|r|/2} \tilde{\gamma}_{J+1}^{(r)}(\sqrt{n} I^{1/2}(\theta_0)(\theta_0 - \hat{\theta}_n))(1 + o_p(1)).$$

Here, $\tilde{\gamma}_{J+1}^{(r)}(\sqrt{n} I^{1/2}(\theta_0)(\theta_0 - \hat{\theta}_n))$ is generated by applying the chain rule to the last term on the right-hand side in (A.1). Note that we differentiate with respect to $\theta$ and then evaluate at $\theta_0$. Expressions (A.3) and (A.4) will give the two leading terms in (3.8), respectively.

Next we use Assumption EE to observe an identity: We can take expectations over $\hat{\theta}_n$ when it occurs in the argument of a polynomial $Q(\cdot)$ by the relationship

$$E_{\theta_0}(Q(\sqrt{n}I^{1/2}(\theta_0)(\theta_0 - \hat{\theta}_n))\phi(\sqrt{n}I^{1/2}(\theta_0)(\theta_0 - \hat{\theta}_n)))$$

(A.6)
$$= \int Q(v)\phi(v)\left(\phi(v) + \sum_{k=1}^{K} n^{-k/2}P_k(v)\phi(v) + \frac{o(n^{-K/2})}{1+\|v\|^{K+2}}\right)dv$$

$$= \frac{1}{(4\pi)^{d/2}}\int Q\left(\frac{v}{\sqrt{2}}\right)\left(\phi(v) + \sum_{k=1}^{K} n^{-k/2}P_k\left(\frac{v}{\sqrt{2}}\right)\phi(v)\right)dv + o(n^{-K/2})$$

$$= \frac{1}{(4\pi)^{d/2}}\left(Q\left(\frac{\sigma}{\sqrt{2}}\right) + \sum_{k=1}^{K} n^{-k/2}Q \circ P_k\left(\frac{\sigma}{\sqrt{2}}\right)\right) + o(n^{-K/2}),$$

where $Q \circ P_k(\cdot)$ is the polynomial obtained by their product in which, as before, we have taken expectations and replaced powers. [The factor $1/(4\pi)^{d/2}$ appears when we multiply two standard normal densities and observe the change of variables in the exponent.]

We use (A.6) in (A.3), (A.4) and (A.5) to get (3.8).

Since the integrability of $W^{(r)}(\cdot|X^n)$ and $H_{r-1}(\cdot)\phi(\cdot)$ implies that of $\tilde{\gamma}_{J+1}(\cdot)$, we can apply (A.6) to see that the expectation of the error term (A.5) is

(A.7)
$$E_{\theta_0}(n^{-(J+1)/2}n^{|r|/2}\tilde{\gamma}_{J+1}^{(r)}(\sqrt{n}I^{1/2}(\theta_0)(\theta_0 - \hat{\theta}_n))(1 + o_p(1)))$$

$$= Cn^{(|r|-J-1)/2}\int \tilde{\gamma}_{J+1}^{(r)}(v)\left(\phi(v) + \sum_{k=1}^{K} n^{-k/2}P_k(v)\phi(v) + o(n^{-K/2})\frac{1}{1+\|v\|^{K+2}}\right)dv$$

$$= O(n^{(|r|-J-1)/2}).$$

In (A.7) we used the fact that the integral over $\tilde{\gamma}_{J+1}^{(r)}(v)\phi(v)$ gives an $O(1)$ term. The integral over the summands in the summation gives terms of order $O(1)n^{k/2}$, for $k=1,\ldots,K$. So, the initial term gives the order in $n$ as indicated in (A.7).

Similarly, using (A.6), the expectation of (A.4) is

$$\sum_{j=1}^{J} n^{(|r|-j)/2}\left[\frac{1+o(1)}{(4\pi)^{d/2}}\left(\eta_j^{(r)}\left(\frac{\sigma}{\sqrt{2}}\right) + \sum_{k=1}^{K} n^{-k/2}\eta_j^{(r)} \circ P_k\left(\frac{\sigma}{\sqrt{2}}\right)\right)\right.$$
$$\left. + o(n^{(|r|-K-1)/2})\right]$$

(A.8)
$$= \sum_{j=1}^{J} n^{(|r|-j)/2} \frac{1+o(1)}{(4\pi)^{d/2}} \eta_j^{(r)}\left(\frac{\sigma}{\sqrt{2}}\right)$$
$$+ \sum_{k+j \leq J} n^{(|r|-k-j)/2} \frac{1+o(1)}{(4\pi)^{d/2}} \eta_j^{(r)} \circ P_k\left(\frac{\sigma}{\sqrt{2}}\right) + o(n^{(|r|-K-1)/2}).$$

The leading term in (A.8) gives the second term in $A_1$ in (3.8).

Finally, using (A.6), the expectation of (A.3) is

(A.9)
$$n^{|r|/2} \frac{|I(\theta_0)|^{1/2}}{(4\pi)^{d/2}} \left( H_{r-1}\left(\frac{\sigma}{\sqrt{2}}\right) \right.$$
$$\left. + \sum_{j=1}^{J} n^{-j/2} H_{r-1} \circ P_j\left(\frac{\sigma}{\sqrt{2}}\right) \right) + o(n^{(|r|-J)/2}),$$

which gives the leading term in (3.8) and the first term in $A_1$. That is, by collecting terms in (A.7)–(A.9), the proof is completed. □

To exemplify Theorem 3.3, we examine the average behavior of the posterior density at $\theta_0$. Straightforward extensions give similar results at other values of $\theta$.

Consider the functional $F(W(\cdot|X^n)) = w(\theta_0|X^n) = \frac{\partial^{|r|} W(\cdot|X^n)}{\partial \theta^r}|_{\theta=\theta_0}$ with $r = (1, \ldots, 1)$. Since $H_{r-1}(\cdot) = H_{\mathbf{0}}(\cdot) \equiv 1$, Theorem 3.3 gives

(A.10)
$$E_{\theta_0}(w(\theta_0|X^n))$$
$$= \frac{n^{d/2} |I(\theta_0)|^{1/2}}{(4\pi)^{d/2}} + n^{(d-1)/2} A_1 + o(n^{(d-1)/2}) + R'_n.$$

When $d = 1$, we can verify that $A_1 = 0$. This is easy because the expressions for the $\gamma_j(\cdot)$'s are available from [11] in this case. Indeed, we have

$$\eta_1(v) = I(\theta_0) c_{00}^{-1}(c_{10} v^3 + c_{01} v)$$

and

$$A_1 = \frac{|I(\theta_0)|^{1/2}}{(4\pi)^{d/2}} P_1\left(\frac{\sigma}{\sqrt{2}}\right) + \frac{1}{(4\pi)^{d/2}} \eta_1^{(r)}\left(\frac{\sigma}{\sqrt{2}}\right),$$

in which $P_1(v) = \chi_3 v/3!$. The expectations of $P_1(v)$ and $\eta_1(v)$ when $v$ is Normal$(0,1)$ are obviously zero. So, $P_1(\frac{\sigma}{\sqrt{2}}) = \eta_1^{(r)}(\frac{\sigma}{\sqrt{2}}) = 0$ and, thus, $A_1 = 0$. This means that the two biggest terms in (A.10) are of order $n^{d/2}$ and $n^{(d-2)/2}$. We have not carried out the analysis far enough to identify the coefficient of the second-order term.

It is seen that (A.10) is the same as the result in [5]. We remark that if one chooses $F(W(\cdot|X^n)) = w^2(\theta_0|X^n)$, the techniques above give

$$(A.11) \quad E_{\theta_0}(w^2(\theta_0|X^n)) \sim E_{\theta_0}(n^d|I(\theta_0)|\phi^2(Z)) = \frac{n^d|(\theta_0)|}{3^{d/2}(2\pi)^d},$$

the same as in [5].

For completeness, we next show how to use the general procedure Proposition 2.1 to get (A.10). There are four types of terms in (2.10); we go through them in turn.

The first term on the right-hand side of (2.10) is

$$EF(\Phi(Z + \sqrt{n}I^{1/2}(\theta_0)(\cdot - \theta_0)))$$
$$= \frac{n^{d/2}|I^{1/2}(\theta_0)|}{(2\pi)^{d/2}} \int \frac{1}{(2\pi)^{d/2}} e^{-(1/2)z'z} e^{-(1/2)z'z}\, dz$$
$$= \frac{n^{d/2}|I^{1/2}(\theta_0)|}{(2\pi)^{d/2}} \int \frac{1}{(2\pi)^{d/2}} e^{-z'z}\, dz = \frac{n^{d/2}|I^{1/2}(\theta_0)|}{(4\pi)^{d/2}}.$$

Next, for $J \geq 1$, the terms in the summation in (2.10) are of the form

$$n^{-j/2} \frac{n^{1/2}|I^{1/2}(\theta_0)|}{(2\pi)^{d/2}} \int \frac{1}{(2\pi)^{1/2}} e^{-(1/2)z'z} P_j(z) e^{-(1/2)z'z}\, dz$$
$$= n^{-j/2} \frac{n^{d/2}|I^{1/2}(\theta_0)|}{(4\pi)^{d/2}} \int e^{-(1/2)z'z} P_j\left(\frac{z}{\sqrt{2}}\right) dz$$
$$= n^{-j/2} \frac{n^{d/2}|I^{1/2}(\theta_0)|}{(4\pi)^{d/2}} P_j\left(\frac{\sigma}{\sqrt{2}}\right),$$

where $P_j(\frac{\sigma}{\sqrt{2}})$ is the expectation of $P_j(\frac{z}{\sqrt{2}})$ with the $z^l$'s replaced by $\sigma_l$'s, the $l$th moments of $N(\mathbf{0}, I_d)$.

Next, for $h(n)$, we observe that

$$\Phi(z + \sqrt{n}I^{1/2}(\theta_0)(\theta - \theta_0))$$
$$= \frac{1}{(2\pi)^{d/2}} \int_{-\infty}^{z+\sqrt{n}I^{1/2}(\theta_0)(\theta-\theta_0)} e^{-(1/2)t't}\, dt$$
$$= \frac{|nI^{1/2}(\theta_0)|^{1/2}}{(2\pi)^{d/2}} \int_{-\infty}^{\theta} e^{-(1/2)(\sqrt{n}I^{1/2}(\theta_0)(v-\theta_0)+z)'(\sqrt{n}I^{1/2}(\theta_0)(v-\theta_0)+z)}\, dv.$$

This gives

$$\left.\frac{\partial^{|r|} \Phi(z + \sqrt{n}I^{1/2}(\theta_0)(\theta - \theta_0))}{\partial \theta^r}\right|_{\theta=\theta_0} = n^{d/2}|I^{1/2}(\theta_0)|\phi(z),$$

and we have

$$h(n) = \int \frac{F(\Phi(z + \sqrt{n}I^{1/2}(\theta_0)(\cdot - \theta_0)))}{1 + \|z\|^J} dz$$

$$= \int \frac{n^{d/2}|I^{1/2}(\theta_0)|\phi(z)}{1 + \|z\|^J} dz$$

$$= \frac{n^{d/2}|I^{1/2}(\theta_0)|}{(2\pi)^{d/2}} \int \frac{e^{-(1/2)z'z}}{1 + \|z\|^J} dz,$$

which is smaller than the leading term when multiplied by $o(n^{J/2})$ for any $J \geq 1$.

It remains to evaluate the expectation of the remainder term. As assumed in the proofs of the theorems, we only need to evaluate it over the "good" sets, and we omit the indicators for them. Write

$$R_n = \frac{d}{d\theta}\left(\sum_{j=1}^{J} n^{-j/2}\phi(\sqrt{n}I^{1/2}(\theta_0)(\theta - \hat{\theta}_n))\right.$$

$$\times \gamma_j(\sqrt{n}I^{-1/2}(\theta_0)(\theta - \hat{\theta}_n))(1 + o(1))$$

$$\left.+ n^{-(J+1)/2}\gamma_{J+1}(\sqrt{n}I^{-1/2}(\theta_0)(\theta - \hat{\theta}_n))(1 + o(1))\right)\bigg|_{\theta=\theta_0}$$

$$= \sum_{j=1}^{J} n^{-j/2}n^{d/2}\phi(\sqrt{n}I^{1/2}(\theta_0)(\theta - \hat{\theta}_n))\eta_j^{(1)}(\sqrt{n}I^{-1/2}(\theta_0)(\theta_0 - \hat{\theta}_n))(1 + o(1))$$

$$+ n^{-(J+1)/2}n^{d/2}|I^{(r)/2}(\theta_0)|\gamma_{J+1}^{(1)}(\sqrt{n}I^{-1/2}(\theta_0)(\theta - \hat{\theta}_n))(1 + o(1)).$$

So, by Assumption EE and (A.6), we get

$$E_{\theta_0} R_n = \sum_{j=1}^{J} \frac{n^{(d-j)/2}}{(4\pi)^{d/2}}\eta_j^{(1)}(\sigma/\sqrt{2})(1 + o(1)) + o(n^{(d-J)/2}),$$

which has lower order than the leading term for $J \geq 1$. Thus, by Proposition 2.1, we get the same result as from Theorem 3.3.

Our final point is that our techniques can be used to approximate the expected value of posterior expectations. Indeed, from (A.1), note that

$$w(\theta|X^n) = \phi_{\hat{\theta}_n, I^{-1}(\theta_0)/n}(\theta)$$

$$+ \sum_{j=1}^{J} n^{-j/2}n^{d/2}|I^{1/2}(\theta_0)|\phi(\sqrt{n}I^{1/2}(\theta_0)(\theta - \hat{\theta}_n))$$

(A.12)

$$\times \eta_j^{(1)}(\sqrt{n}I^{1/2}(\theta_0)(\theta - \hat{\theta}_n))(1 + o(1))$$

$$+ n^{-(J+1)/2}n^{d/2}|I^{1/2}(\theta_0)|\gamma_{J+1}^{(1)}(\sqrt{n}I^{1/2}(\theta_0)(\theta - \hat{\theta}_n))(1+o(1)).$$

The $\gamma_j(\cdot, X^n)$'s are from Assumption JE and are differentiable, as are the $\eta_j^{(r)}(\cdot)$'s. Now, suppose $h = h(\theta)$ has all $r$th partial derivatives, for $|r| \leq J$, on a neighborhood of $\theta_0$ and that $h(\theta)w(\theta|X^n)$ and its partial derivatives are integrable with respect to $w(\cdot|X^n)$.

Then, Taylor expanding $h$ at $\hat{\theta}$, justifying a use of Assumption EE and gathering terms suggests that

$$E_{\theta_0}\left(\int h(\theta)w(\theta|X^n)\,d\theta\right) = h(\theta_0) + n^{-1/2}I^{-1/2}(\theta_0)\sum_{|r|=1}h^{(r)}(\theta_0)\sigma_r$$
(A.13)
$$+ A_1 n^{-1} + o(n^{-1}) + R'_n,$$

where

(A.14) $\quad A_1 = I^{-1/2}(\theta_0)\sum_{|r|=1}\eta_1^{(r)}(\sigma)h^{(r)}(\theta_0) + \frac{3}{2}I^{-1}(\theta_0)\sum_{|r|=2}h^{(r)}(\theta_0)$

and the $\eta_j^{(r)}(\cdot)$'s are as in Theorem 3.3. An extension of this argument gives similar expressions for higher-order terms.

**Acknowledgments.** Most of this work was done while the first author was on leave at ISDS, Duke University and SAMSI in Research Triangle Park.

DEPARTMENT OF STATISTICS
UNIVERSITY OF BRITISH COLUMBIA
333-6356 AGRICULTURAL ROAD
VANCOUVER, BRITISH COLUMBIA
CANADA V6T 1Z2
E-MAIL: riffraff@stat.ubc.ca

NATIONAL HUMAN GENOME CENTER
HOWARD UNIVERSITY
2216 SIXTH STREET, N.W., SUITE 205
WASHINGTON, DISTRICT OF COLUMBIA 20059
USA
E-MAIL: ayuan@howard.edu